\documentclass[10pt,a4paper]{article}
\usepackage[a4paper]{geometry}
\usepackage{amssymb,latexsym,amsmath,amsfonts,amsthm}
\usepackage{graphicx}

\usepackage{epsfig}

\DeclareMathOperator{\diag}{diag} 

\newcommand{\er}{\mathbb{R}}
\newcommand{\cee}{\mathbb{C}}
\newcommand{\enn}{\mathbb{N}}

\newcommand{\bewijs}{\textsc{proof}}
\newcommand{\bol}{\hfill\square\\}

\newcommand{\Pee}{{\cal P}}

\newcommand{\vece}{\mathbf{e}}

\renewcommand{\vec}{\mathbf}

\newcommand{\vecm}{\mathbf{m}}
\newcommand{\vecn}{\mathbf{n}}
\newcommand{\vecw}{\mathbf{w}}
\newcommand{\vecA}{\mathbf{A}}
\newcommand{\vecB}{\mathbf{B}}
\newcommand{\vecP}{\mathbf{P}}
\newcommand{\vecQ}{\mathbf{Q}}

\newtheorem{theorem}{Theorem}[section]
\newtheorem{lemma}[theorem]{Lemma}
\newtheorem{proposition}[theorem]{Proposition}

\newtheorem{corollary}[theorem]{Corollary}

\newtheorem{rhp}[theorem]{RH problem}

\theoremstyle{definition}
\newtheorem{definition}[theorem]{Definition}

\theoremstyle{remark}

\newtheorem{remark}[theorem]{Remark}

\numberwithin{equation}{section}

\title{Average characteristic polynomials for multiple orthogonal polynomial ensembles}

\author{Steven Delvaux\footnotemark[1]}
\date{}

\begin{document}

\maketitle
\renewcommand{\thefootnote}{\fnsymbol{footnote}}
\footnotetext[1]{Department of Mathematics, Katholieke Universiteit Leuven,
Celestijnenlaan 200B, B-3001 Leuven, Belgium. email:
steven.delvaux\symbol{'100}wis.kuleuven.be. The author is a Postdoctoral Fellow
of the Fund for Scientific Research - Flanders (Belgium).}

\begin{abstract}
Multiple orthogonal polynomials (MOP) are a non-definite version of matrix
orthogonal polynomials. They are described by a Riemann-Hilbert matrix $Y$
consisting of four blocks $Y_{1,1}$, $Y_{1,2}$, $Y_{2,1}$ and $Y_{2,2}$. In
this paper, we show that $\det Y_{1,1}$ ($\det Y_{2,2}$) equals the average
characteristic polynomial (average inverse characteristic polynomial,
respectively) over the probabilistic ensemble that is associated to the MOP. In
this way we generalize classical results for orthogonal polynomials, and also
some recent results for MOP of type~I and type~II. We then extend our results
to arbitrary products and ratios of characteristic polynomials. In the latter
case an important role is played by a matrix-valued version of the
Christoffel-Darboux kernel. Our proofs use determinantal identities involving
Schur complements, and adaptations of the classical results by Heine,
Christoffel and Uvarov.


\end{abstract}

\section{Introduction and statement of results}
\label{section:introduction}

\subsection{Random matrix ensembles}
\label{subsection:randommxintro}

On the space $\mathcal H_n$ of Hermitian $n$ by $n$ matrices consider the
random matrix ensemble defined by the probability distribution
\begin{equation}\label{randommxensemble}
\frac{1}{Z_n} e^{-\textrm{Tr} V(M)}\ dM,\qquad M\in\mathcal H_n,
\end{equation}
for some given polynomial $V$ of even degree. Here $Z_n$ is a normalization
constant, $\textrm{Tr}$ denotes the trace and $dM$ is the Lebesgue measure on
$\mathcal H_n$.

The random matrix ensemble \eqref{randommxensemble} leads to a probability
distribution on the space $$\er^n_{\leq} = \{(x_1,\ldots,x_n)\in\er^n|\ x_1\leq
\ldots\leq x_n\}$$ of ordered eigenvalue tuples. It is well-known that this
joint probability distribution has the form
\begin{equation}\label{jointPDF1} \frac{1}{\tilde Z_n}\prod_{i<j} (x_j-x_i)^2 \prod_{j=1}^{n}
e^{-V(x_j)} \prod_{j=1}^{n} dx_j,
\end{equation}
for some normalization constant $\tilde Z_n$. Thus the probability that the
ordered eigenvalues of the matrix \eqref{randommxensemble} lie in an
infinitesimal box
$[x_1,x_1+dx_1]\times\ldots\times[x_n,x_n+dx_n]\subset\er_{\leq}^n$ is given by
\eqref{jointPDF1}. Note in particular that the density \eqref{jointPDF1} is
small if two eigenvalues $x_i$ and $x_j$ are close to each other. This means
that the eigenvalues tend to \lq repel\rq\ each other.

One can write \eqref{jointPDF1} alternatively as
\begin{equation}\label{jointPDF2} \frac{1}{\tilde Z_n}\det \left( f_i(x_j)
\right)_{i,j=1}^n\det \left( g_i(x_j) \right)_{i,j=1}^n\ \prod_{j=1}^{n} dx_j
\end{equation}
where \begin{equation}\label{jointPDF2bis} f_i(x) = x^{i-1},\qquad g_i(x) =
x^{i-1}e^{-V(x)}.\end{equation} Indeed, this follows upon recognizing
\eqref{jointPDF1} to be basically a product of two Vandermonde determinants.

Rather than the space $\er^n_{\leq}$ of ordered eigenvalue tuples, we will find
it convenient to consider the probability distribution
\eqref{jointPDF2}--\eqref{jointPDF2bis} on the full space $\er^n$. To maintain
a probability distribution we should then multiply the normalization constant
$\tilde Z_n$ by a factor $n!$.

To the probability distribution \eqref{jointPDF2}--\eqref{jointPDF2bis} one can
associate the \emph{average characteristic polynomial}
\begin{equation}\label{averagecharpol:intro} P_n(z) = \frac{1}{\tilde
Z_{n}}\int_{-\infty}^{\infty}\ldots\int_{-\infty}^{\infty} \left(\prod_{j=1}^n
(z-x_j)\right)\ 
\det \left( f_i(x_j) \right)_{i,j=1}^n\det \left( g_i(x_j) \right)_{i,j=1}^n
\prod_{j=1}^{n} dx_j.
\end{equation}
It follows from a classical calculation of Heine, see e.g.\ \cite{Dei,Sz}, that
$P_n$ can be characterized as the $n$th monic \emph{orthogonal polynomial} with
respect to the weight function $e^{-V(x)}$ on~$\er$. Thus the polynomials $P_n$
satisfy the conditions
$$ P_n(x) = x^n+O(x^{n-1})
$$
for all $n\in\enn$ and \begin{equation}\label{def:cn} \int_{-\infty}^{\infty}
P_n(x)P_m(x)e^{-V(x)}\ dx = c_nc_m\delta_{m,n}
\end{equation}
for all $n,m\in\enn$, for certain $c_n\in\er$.

Intuitively, the above result states that the zeros of the monic orthogonal
polynomial $P_n(z)$ determine the \lq typical\rq\ eigenvalue configuration of
the random matrix ensemble \eqref{randommxensemble}. This holds in particular
in the Gaussian case $V(x) = x^2/2$, where the ensemble
\eqref{randommxensemble} reduces to the well-known \emph{Gaussian unitary
ensemble} (GUE) while the corresponding orthogonal polynomials are (up to
scaling factors) the classical \emph{Hermite polynomials}.

In the literature, many more results on average characteristic polynomials for
random matrix ensembles can be found, see e.g.\ \cite{BDS,BS,BH3,KS,SF} and the
references therein. We mention the following result of Fyodorov-Strahov
\cite{FS}. Define the \emph{average inverse characteristic polynomial}
corresponding to the probability distribution
\eqref{jointPDF2}--\eqref{jointPDF2bis} as
\begin{equation}\label{averageinvcharpol:intro} Q_n(z) = \frac{1}{\tilde
Z_{n}}\int_{-\infty}^{\infty}\ldots\int_{-\infty}^{\infty} \left(\prod_{j=1}^n
(z-x_j)^{-1}\right)\ 
\det \left( f_i(x_j) \right)_{i,j=1}^n\det \left( g_i(x_j) \right)_{i,j=1}^n
\prod_{j=1}^{n} dx_j,
\end{equation}
for $z\in\cee\setminus\er$. Then it holds that
\begin{equation}\label{Strahov:intro} Q_n(z) = \frac{1}{c^2_{n-1}}\int_{-\infty}^{\infty} \frac{P_{n-1}(x)e^{-V(x)}}{z-x}\
dx,
\end{equation}
with $c_n$ as in \eqref{def:cn}. Fyodorov--Strahov \cite{FS,SF} also observed
that both $P_n(z)$ and the right hand side of \eqref{Strahov:intro} have a
natural interpretation in terms of the associated \emph{Riemann-Hilbert
problem} (briefly \emph{RH problem}), and they obtained determinantal formulae
for averages of more general products and ratios of characteristic polynomials;
see further.

In recent years there has been interest in the following generalization of
\eqref{randommxensemble},
\begin{equation}\label{randommx:source}
\frac{1}{Z_n} e^{-\textrm{Tr}(V(M)-AM)}\ dM.
\end{equation}
Here $A$ is a fixed diagonal matrix which is called the \emph{external source}.
Typically $A$ has only a small number of distinct eigenvalues $a_1,\ldots,a_p$,
say with corresponding multiplicities $n_1,\ldots,n_p$.

The model \eqref{randommx:source} was first studied by Br\'ezin-Hikami
\cite{BH1,BH2} and P. Zinn-Justin \cite{Zinn} who showed that the eigenvalue
correlations are determinantal. In \cite{BK1} it was shown that
the eigenvalues $x_1,\ldots,x_n$ have a joint probability distribution
of the form \eqref{jointPDF2}, where now $f_i(x) = x^{i-1}$ (Vandermonde
factor) while $g_i(x)$ are given by several Vandermonde-like series:
\begin{align*}
g_i(x) & = x^{i-1} e^{-(V(x)-a_1 x)}, \qquad i = 1, \ldots, n_1, \\
g_{n_1+i}(x) & = x^{i-1} e^{-(V(x)-a_2 x)}, \qquad i = 1, \ldots, n_2, \\
\vdots & \\
g_{n_1+\ldots+n_{p-1}+i}(x) & = x^{i-1}e^{-(V(x)-a_p x)}, \qquad i = 1, \ldots,
n_p.
\end{align*}
To this eigenvalue ensemble one can associate the average characteristic
polynomial $P_n(z)$ in exactly the same way as before, see
\eqref{averagecharpol:intro}. Bleher-Kuijlaars \cite{BK1} showed that $P_n(z)$
satisfies the multiple orthogonality relations
$$ \int_{-\infty}^{\infty} P_n(x)x^{i} e^{-(V(x)-a_k x)}\ dx = 0,
$$
for $i=0,\ldots,n_k-1$ and $k=1,\ldots,p$. For more background on this kind of
orthogonality relations, see e.g.\ \cite{Apt,VA}. Desrosiers-Forrester
\cite{Desrosiers1} showed that the result \eqref{Strahov:intro} on the average
inverse characteristic polynomial \eqref{averageinvcharpol:intro} can also be
generalized to the random matrix ensemble with external source
\eqref{randommx:source}.

The goal of this paper is to generalize the above results to the more general
context of \emph{multiple orthogonal polynomial ensembles} (briefly \emph{MOP
ensembles}) in the sense of Daems-Kuijlaars \cite{DK2}. We will show that in
general, \eqref{averagecharpol:intro} and \eqref{averageinvcharpol:intro} can
be expressed as \emph{Riemann-Hilbert minors}, i.e., determinants of certain
submatrices of the Riemann-Hilbert matrix. Using these results, we will then
obtain determinantal formulas for averages of arbitrary products and ratios of
characteristic polynomials, by adapting the method of Baik-Deift-Strahov
\cite{BDS}. Our results will be stated in Section \ref{subsection:mainresults}.
In Sections \ref{subsection:detpointintro}--\ref{subsection:RHintro} we first
recall the basic definitions concerning MOP ensembles.

\subsection{MOP ensembles}
\label{subsection:detpointintro}

We consider a stochastic model in the following way \cite{DK2}. Let
$p,q\in\mathbb N $ be two positive integers. Let there be given
\begin{itemize}
\item A (finite) sequence of positive integers $n_1,n_2,\ldots,n_p\in\mathbb N $;
\item A sequence of weight functions $w_{1,1}(x),w_{1,2}(x),\ldots,w_{1,p}(x):\mathbb R \to\mathbb R $;
\item A sequence of positive integers $m_1,m_2,\ldots,m_q\in\mathbb N $;
\item A sequence of weight functions
$w_{2,1}(x),w_{2,2}(x),\ldots,w_{2,q}(x):\mathbb R \to\mathbb R $.
\end{itemize}
We will use the vector notations $\vecn := (n_1, \ldots, n_p)$, $|\vecn| :=
\sum_{k=1}^p n_k$ and similarly for $\vecm$ and $|\vecm|$. Occasionally we will
also write $\vecw_1(x) := (w_{1,1}(x),\ldots,w_{1,p}(x))$ and similarly for
$\vecw_2(x)$. Assume in what follows that $|\vecn| = |\vecm| =: n$.

To the above data one can associate a stochastic model, called a \emph{multiple
orthogonal polynomial ensemble} or briefly \emph{MOP ensemble}. This model
consists of $n$ random points $x_1,\ldots,x_n$ on the real line whose joint
p.d.f.\ can be written as a product of two determinants:

\begin{equation} \label{determinantalpointprocess}
    \frac{1}{Z_{n}} \det \left( f_i(x_j) \right)_{i,j=1}^n
        \cdot \det \left( g_i(x_j) \right)_{i,j=1}^n,
\end{equation}
where $Z_{n}$ is a normalization factor, and with functions
\begin{align}
    \nonumber f_i(x) & = x^{i-1} w_{1,1}(x), \qquad i = 1, \ldots, n_1, \\
    \nonumber f_{n_1+i}(x) & = x^{i-1} w_{1,2}(x), \qquad i = 1, \ldots, n_2, \\
    \label{def:fi} \vdots & \\
    \nonumber f_{n_1+\ldots+n_{p-1}+i}(x) & = x^{i-1} w_{1,p}(x), \qquad i = 1, \ldots, n_p,
\end{align}
and
\begin{align}
    \nonumber g_i(x) & = x^{i-1} w_{2,1}(x), \qquad i = 1, \ldots, m_1, \\
    \nonumber g_{m_1+i}(x) & = x^{i-1} w_{2,2}(x), \qquad i = 1, \ldots, m_2, \\
    \label{def:gi} \vdots & \\
    \nonumber g_{m_1+\ldots+m_{q-1}+i}(x) & = x^{i-1} w_{2,q}(x), \qquad i = 1, \ldots,
    m_q.
\end{align}
Note that this is a special case of a biorthogonal ensemble \cite{Bor}. The
motivation for the name \lq MOP ensemble\rq\ comes from the multiple orthogonal
polynomials that are related to it, see Section \ref{section:MOP}.

In order for the above model to be a valid probability distribution on $\er^n$
one should have that \eqref{determinantalpointprocess} is positive on $\er^n$.
We will not address this topic here since our algebraic results will be valid
irrespective of this positivity condition.

The normalization constant $Z_n$ in \eqref{determinantalpointprocess} serves to
make the total probability of the MOP ensemble on $\er^n$ equal to $1$. We will
further obtain an expression for $Z_n$ as a block Hankel determinant; see the
remark at the end of Section \ref{subsection:prooftheorem:ratio}.

\begin{figure}[t]
\begin{center}\vspace{-1mm}
\includegraphics[scale=0.65]{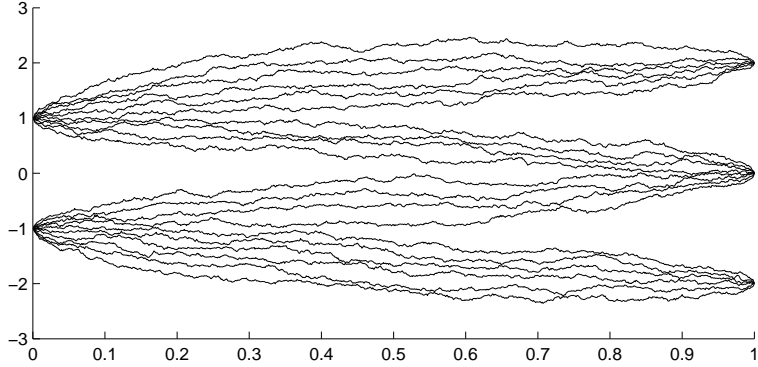}
\end{center}\vspace{-5mm}
\caption{The figure shows $n=20$ non-intersecting Brownian motions with $p=2$
starting points and $q=3$ ending points. The horizontal axis denotes the time
$t\in[0,1]$ and the vertical axis denotes the space variable $x$. For each
$t\in(0,1)$ the positions $x_1,\ldots,x_n$ of the paths at time $t$ form a MOP
ensemble \eqref{determinantalpointprocess}--\eqref{def:gi}.} \label{figwigner4}
\end{figure}

Of course, the above model contains the eigenvalue ensembles in Section
\ref{subsection:randommxintro} as special cases. Further motivation comes from
the theory of non-intersecting one-dimensional Brownian paths with $p$ distinct
starting positions and $q$ distinct ending positions, see Figure
\ref{figwigner4}. See \cite{DK2,KMcG} for details. Another type of application
can be found in \cite{Fid}.

To facilitate comparison with the literature, we may note the following
terminology which is often used: the MOP ensemble is said to be of \emph{type
I} if $q=1$, of \emph{type II} if $p=1$ and of \emph{mixed type} if $p,q\geq
2$. Thus the results of our paper will be valid for general MOP ensembles of
mixed type, whereas the results of \cite{BK1,Desrosiers1,Kui} only apply to MOP
ensembles of type I or type II.

\subsection{Riemann-Hilbert problem and kernel}
\label{subsection:RHintro}

The MOP ensemble in Section \ref{subsection:detpointintro} and its associated
multiple orthogonal polynomials are intimately related to the following
Riemann-Hilbert problem (RH problem) introduced in \cite{DK2}. The RH problem
generalizes the well-known RH problem for orthogonal polynomials due to
Fokas-Its-Kitaev \cite{FIK} as well as its generalization in~\cite{VAGK}.

\begin{rhp}\label{RHP:original}
Consider the data $\vecn$, $\vecw_1(x)$, $\vecm$ and $\vecw_2(x)$ as above and
assume that $|\vecn|=|\vecm|=:n$. The RH problem consists in finding a
matrix-valued function $Y(z)=Y_{\vecn,\vecw_1,\vecm,\vecw_2}(z)$ of size $p+q$
by $p+q$ such that
\begin{itemize}
\item[(1)] $Y(z)$ is analytic
in $\mathbb C \setminus\mathbb R $;
 \item[(2)] For $x\in\mathbb R $, it holds that
\begin{equation}\label{jumpsY10}
Y_{+}(x) = Y_{-}(x)
\begin{pmatrix} I_p & W(x)\\
0 & I_q
\end{pmatrix},
\end{equation}
where $I_k$ denotes the identity matrix of size $k$; where
\begin{equation}\label{weightmatrixrankone}
W(x) = \begin{pmatrix}w_{1,1}(x)w_{2,1}(x)&\ldots&w_{1,1}(x)w_{2,q}(x)\\
\vdots& & \vdots \\
w_{1,p}(x)w_{2,1}(x)&\ldots&w_{1,p}(x)w_{2,q}(x)\end{pmatrix} =
\begin{pmatrix} w_{1,1}(x)\\ \vdots \\ w_{1,p}(x) \end{pmatrix}
\begin{pmatrix} w_{2,1}(x)& \ldots & w_{2,q}(x) \end{pmatrix},
\end{equation}
and where the notation $Y_+(x), Y_-(x)$ denotes the limit of $Y(z)$ with $z$
approaching $x\in\mathbb R $ from the upper or lower half plane in $\mathbb C
$, respectively;
  \item[(3)] As $z\to\infty$, we have that
\begin{equation}\label{asymptoticconditionY0} Y(z) =
    (I_{p+q}+O(1/z))\diag(z^{n_1},\ldots,z^{n_p},z^{-m_1},\ldots,z^{-m_q}).
\end{equation}
\end{itemize}
\end{rhp}

The solution $Y(z)$ to the RH problem is unique, if it exists. Partition this
matrix as
\begin{equation}\label{partitionY} Y(z) =
\begin{array}{ll} & \begin{array}{ll} \ \ \ \ p & \ \ \ \ \ \ \ \ q \end{array}\\
\begin{array}{l} p \\ q\end{array} \hspace{-4mm} & \begin{pmatrix}Y_{1,1}(z) & Y_{1,2}(z)\\ Y_{2,1}(z) &
Y_{2,2}(z)
\end{pmatrix},\end{array}
\end{equation}
where the partition is such that $Y_{1,1}(z)$ has size $p\times p$, and so on.
Then the entries of $Y_{1,1}$ and $Y_{2,1}$ can be described in terms of
multiple orthogonal polynomials, while $Y_{1,2}$ and $Y_{2,2}$ contain certain
Cauchy transforms thereof; see Section \ref{section:MOP} for details.

The RH problem can be used to define the \emph{Christoffel-Darboux kernel}.
There are actually two kernels into play. The first is a $q\times p$
matrix-valued kernel $K_n(x,y)$ defined as
\begin{equation}\label{def:kernel:block}
K_n(x,y) = \frac{1}{2\pi i(y-x)}\begin{pmatrix}0 & I_q \end{pmatrix}
Y^{-1}(x)Y(y)
\begin{pmatrix}I_p \\ 0 \end{pmatrix},
\end{equation}
for $x,y\in\cee$, $x\neq y$. This definition may seem rather unmotivated; but
see Section \ref{section:MOP} for the connection with the classical definition
of the Christoffel-Darboux kernel.

Second, there is also a scalar-valued kernel $\hat K_n(x,y)$ defined by
\begin{equation}\label{def:kernel:scalar} \hat K_n(x,y)
=
\begin{pmatrix}w_{2,1}(x) & \ldots & w_{2,q}(x)\end{pmatrix} K_n(x,y)
\begin{pmatrix}w_{1,1}(y) \\ \vdots \\ w_{1,p}(y)\end{pmatrix}.
\end{equation}
It is known that the MOP ensemble in Section \ref{subsection:detpointintro} is
\emph{determinantal} with \emph{correlation kernel} $\hat{K}_n(x,y)$
\cite{DK2}. In this paper, however, the quantity of interest will be the
matrix-valued kernel $K_n(x,y)$.


\subsection{Statement of results}
\label{subsection:mainresults}

Now we are ready to state our main results. Throughout this section we assume
fixed data $\vecn$, $\vecw_1(x)$, $\vecm$ and $\vecw_2(x)$ as before, and we
consider the corresponding MOP ensemble
\eqref{determinantalpointprocess}--\eqref{def:gi}. For any $K,L\in\enn$ we
define $P_{n}^{[K,L]}$ as
\begin{multline}\label{averageproductratioKL} P_{n}^{[K,L]}(y_1,\ldots,y_K;z_1,\ldots,z_L) :=
\\
\frac{1}{Z_{n}}\int_{-\infty}^{\infty}\ldots\int_{-\infty}^{\infty}
\frac{\prod_{k=1}^K\left(\prod_{j=1}^n
(y_k-x_j)\right)}{\prod_{k=1}^L\left(\prod_{j=1}^n (z_k-x_j)\right)} \det
\left( f_i(x_j) \right)_{i,j=1}^n
        \cdot \det \left( g_i(x_j) \right)_{i,j=1}^n\ dx_1\ldots dx_n,
\end{multline}
i.e., $P_{n}^{[K,L]}$ is the average with respect to the MOP ensemble
\eqref{determinantalpointprocess}--\eqref{def:gi} of the ratio of products of
characteristic polynomials. Here we assume that $y_1,\ldots,y_K\in\cee$,
$z_1,\ldots,z_L\in\cee\setminus\er$ and the numbers in the set
$(y_1,\ldots,y_K$, $z_1,\ldots,z_L)$ are pairwise distinct.

Note that for $K=1$ and $L=0$, \eqref{averageproductratioKL} reduces to the
average characteristic polynomial \eqref{averagecharpol:intro}, while for $K=0$
and $L=1$ it reduces to the average inverse characteristic polynomial
\eqref{averageinvcharpol:intro}. These are the cases under consideration in our
first two main theorems.

\begin{theorem}\label{theorem:main1} (Average characteristic polynomial)
We have that
\begin{equation}\label{mainP3} P^{[1,0]}_n(y) = \det Y_{1,1}(y),
\end{equation}
according to the partition \eqref{partitionY}.
\end{theorem}

\begin{theorem}\label{theorem:main2} 
(Average inverse characteristic polynomial) We have that
\begin{equation}\label{mainQ3}
P_n^{[0,1]}(z) = \det Y_{2,2}(z),
\end{equation}
according to the partition \eqref{partitionY}.
\end{theorem}

The proofs of Theorems \ref{theorem:main1}--\ref{theorem:main2} will be given
in Section~\ref{subsection:prooftheorem:1and2}.

\begin{remark}\label{remark:symmetry}
\begin{itemize}\textrm{}
\item[(a)]
One may understand Theorem \ref{theorem:main1} as follows. The quantity $\det
Y_{1,1}(y)$ is a monic polynomial of degree $n$, and its zeros determine the
\lq typical\rq\ point configuration of the MOP ensemble in Section
\ref{subsection:detpointintro}. By specializing this to the MOP ensembles in
the introduction, one obtains the typical eigenvalue configuration of random
matrices with external source, the typical positions of non-intersecting
Brownian motions, etc. However, one should be careful with these
interpretations. For one thing, we do not even now if all the zeros of $\det
Y_{1,1}(y)$ are real. This is an open problem.
\item[(b)] The functions $P_n^{[K,L]}$ in \eqref{averageproductratioKL} are clearly invariant under the involution
$(\vecn,\vecw_1)\leftrightarrow (\vecm,\vecw_2)$. The similar fact for the
right hand sides of \eqref{mainP3} and \eqref{mainQ3} is discussed in
Section~\ref{section:MOP}, see \eqref{jacobi1}--\eqref{jacobi2}.
\end{itemize}
\end{remark}

As we mentioned before, Theorem \ref{theorem:main1} generalizes a classical
result for orthogonal polynomials, see e.g.\ \cite{Dei}, as well as its
generalization to type II MOP due to
Bleher-Kuijlaars \cite{BK1}. 
On the other hand, Theorem \ref{theorem:main2} generalizes a result for
orthogonal polynomials due to Fyodorov-Strahov \cite{FS}, and also its
generalization to type I MOP due to Desrosiers-Forrester \cite{Desrosiers1},
see also \cite{Kui}.

Now we state our third main theorem.

\begin{theorem}\label{theorem:main3} (Average ratio of characteristic polynomials)
For the case where $K=L=1$ in \eqref{averageproductratioKL}, we have that
\begin{equation}\label{mainR3b}
P_n^{[1,1]}(y,z) = \det L_n(y,z)
\end{equation}
where $L_n$ is the $q$ by $q$ matrix
\begin{multline}\label{def:Ln} L_n(y,z) := I_q-(z-y)\int_{-\infty}^{\infty} K_n(y,x)W(x)/(z-x)\ dx
= \begin{pmatrix}0 & I_q
\end{pmatrix}Y^{-1}(y)Y(z)\begin{pmatrix}0 \\ I_q
\end{pmatrix}.
\end{multline}
Alternatively, we also have
\begin{equation}\label{mainR3a}
P_n^{[1,1]}(y,z) = \det R_n(z,y)
\end{equation}
where $R_n$ is the $p$ by $p$ matrix
\begin{multline}\label{def:Rn} R_n(z,y) := I_p-(z-y)\int_{-\infty}^{\infty} W(x)K_n(x,y)/(z-x)\ dx
= \begin{pmatrix}I_p & 0 \end{pmatrix}Y^{-1}(z)Y(y)\begin{pmatrix}I_p \\
0
\end{pmatrix}.
\end{multline}
Here $W(x)$ is the weight matrix \eqref{weightmatrixrankone} and $K_n(x,y)$ is
the Christoffel-Darboux kernel \eqref{def:kernel:block}.
\end{theorem}

The proof of Theorem \ref{theorem:main3} will be given in Section
\ref{subsection:prooftheorem:ratio}.

Note that the matrices $W(x)$ and $K_n(x,y)$ are of size $p\times q$ and
$q\times p$, respectively, so each of the matrix expressions in \eqref{def:Ln}
and \eqref{def:Rn} is well-defined. (The integrals are taken entry-wise.) The
equalities between these expressions are discussed in Section
\ref{subsection:kernel}.

\begin{corollary}\label{corollary:main3} Assume that $q=1$ and $w_{2,1}(x)\equiv 1$.
Then
\begin{equation}\label{corollary:main3a} P_n^{[1,1]}(y,z) = 1-(z-y)\int_{-\infty}^{\infty} \hat K_n(y,x)/(z-x)\
dx,
\end{equation} where $\hat K_n$ is the scalar kernel in
\eqref{def:kernel:scalar}. In particular,
\begin{equation}\label{corollary:main3b}
\hat K_n(y,z) = \frac{1}{z-y}\lim_{\epsilon\to 0} \frac{1}{2\pi i}
\left(P_n^{[1,1]}(y,z+\epsilon i)-P_n^{[1,1]}(y,z-\epsilon i)\right),\qquad
z\in\er.
\end{equation}
Similar results hold when $p=1$ and $w_{1,1}(x)\equiv 1$.
\end{corollary}

Formula \eqref{corollary:main3b} 
retrieves a result in \cite{Desrosiers1}, see also \cite{BS} for the scalar
case $p=q=1$.


Finally, let us consider the case of arbitrary $K$ and $L$ in
\eqref{averageproductratioKL}. To this end we need some extra notation.

\begin{definition}\label{def:chainmultiindices} Assume that $|\vecn| = |\vecm|$. We
set $(\vecn_0,\vecm_0):=(\vecn,\vecm)$ and we define the sequence of
multi-indices $(\vecn_k,\vecm_k)$, each $\vecn_k$ being a vector of length $p$
and each $\vecm_k$ a vector of length $q$, for $k\in\mathbb Z$ inductively as
follows
\begin{itemize}
\item For $k=1,2,\ldots$, we set $\vecn_{k} := \vecn_{k-1}+ (1,1,\ldots,1)$
and we fix $\vecm_{k}$ arbitrarily such that $\vecm_{k}\geq \vecm_{k-1}$
componentwise and $|\vecm_{k}-\vecm_{k-1}|=p$.
\item For
$k=-1,-2,\ldots$, we set $\vecm_{k} := \vecm_{k+1}- (1,1,\ldots,1)$ and we fix
$\vecn_{k}$ arbitrarily such that $\vecn_{k}\leq \vecn_{k+1}$ componentwise and
$|\vecn_{k+1}-\vecn_{k}|=q$.
\end{itemize}
\end{definition}

Definition \ref{def:chainmultiindices} implies that $|\vecn_k| = |\vecm_k|$ for
each $k$, so we can consider the RH problem with respect to the pair of
multi-indices $(\vecn_k,\vecm_k)$. We will assume that the RH problem is
solvable for all involved pairs of multi-indices.

Note that in Definition \ref{def:chainmultiindices} the multi-indices
$(\vecn_k,\vecm_k)$ for negative $k$ are only meaningful provided $|k|\leq
\min_{l=1,\ldots,q} m_l$. Otherwise we would obtain multi-indices with negative
components.

Now we are ready to consider \eqref{averageproductratioKL} for arbitrary $K$
and $L$. It turns out that the results can be expressed as determinants of
large matrices constructed from the \lq building blocks\rq\ in Theorems
\ref{theorem:main1}--\ref{theorem:main3}. To get the idea we first describe the
three typical cases.

\begin{theorem}\label{theorem:main4} For the expression $P_n^{[K,L]}$ in
\eqref{averageproductratioKL}, we have
\begin{itemize}
\item[(a)]
If $L=0$ then
\begin{equation}
P_{n}^{[K,0]}(y_1,\ldots,y_K) = \frac{1}{\prod_{1\leq i<j\leq K}
(y_j-y_i)^p}\det\begin{pmatrix}
(Y_{1,1})_{\vecn_0,\vecm_0}(y_1) & \ldots & (Y_{1,1})_{\vecn_0,\vecm_0}(y_K) \\
\vdots & & \vdots \\
(Y_{1,1})_{\vecn_{K-1},\vecm_{K-1}}(y_1) & \ldots &
(Y_{1,1})_{\vecn_{K-1},\vecm_{K-1}}(y_K)
\end{pmatrix}.
\end{equation}
\item[(b)] If $K=0$ then
\begin{equation}
P_{n}^{[0,L]}(z_1,\ldots,z_L) = \frac{1}{\prod_{1\leq i<j\leq L}
(z_j-z_i)^q}\det\begin{pmatrix}
(Y_{2,2})_{\vecn_0,\vecm_0}(z_1) & \ldots & (Y_{2,2})_{\vecn_0,\vecm_0}(z_L) \\
\vdots & & \vdots \\
(Y_{2,2})_{\vecn_{-(L-1)},\vecm_{-(L-1)}}(z_1) & \ldots &
(Y_{2,2})_{\vecn_{-(L-1)},\vecm_{-(L-1)}}(z_L)
\end{pmatrix}.
\end{equation}
Here we assume that $L\leq 1+\min_{l=1}^{q} m_l$.
\item[(c)] If $K=L$ then
\begin{multline}\label{main4c}
P_{n}^{[K,K]}(y_1,\ldots,y_K,z_1,\ldots,z_K) \\ = \left(\frac{\prod_{i,j=1}^{K}
(z_j-y_i)}{\prod_{1\leq i<j\leq K} (y_j-y_i)(z_i-z_j)}\right)^p
\det\begin{pmatrix}
\frac{1}{z_1-y_1} R_n(z_1,y_1) & \ldots & \frac{1}{z_1-y_K} R_n(z_1,y_K) \\
\vdots & & \vdots \\
\frac{1}{z_K-y_1} R_n(z_K,y_1) & \ldots & \frac{1}{z_K-y_K} R_n(z_K,y_K)
\end{pmatrix},
\end{multline}
where $R_n$ is the matrix in \eqref{def:Rn}. A similar result holds with the
matrix $L_n$ in \eqref{def:Ln}.
\end{itemize}
\end{theorem}

Theorem \ref{theorem:main4} generalizes results for the scalar case in
\cite{FS}. Our proof will be based on the methods of \cite{BDS}.

We note that Parts (a)--(b) of Theorem \ref{theorem:main4} are of a different
flavor than Part (c). In fact, it is shown in \cite{BDS,FS} that in the scalar
case, 
one can obtain \lq mixtures\rq\ of the formulas in Parts (a)--(b); but these
formulas do not seem to have an analogue in the present setting. This is the
reason why we need to work with the formula of \lq two-point\rq\ type in Part
(c).

The case of arbitrary $K,L$ in \eqref{averageproductratioKL} is obtained from
mixtures of either Parts (a) and (c), or Parts (b) and (c) of Theorem
\ref{theorem:main4}.

\begin{theorem}\label{theorem:main5} For the expression $P_n^{[K,L]}$ in
\eqref{averageproductratioKL}, we have \begin{itemize}
\item[(a)] If $K\geq L$ then
\begin{multline}
P_{n}^{[K,L]}(y_1,\ldots,y_K,z_1,\ldots,z_L) =
\left(\frac{(-1)^{L(K-L)}\prod_{i,j}
(z_j-y_i)}{\prod_{1\leq i<j\leq K} (y_j-y_i)\prod_{1\leq i<j\leq L} (z_i-z_j)}\right)^p\times \\
\det\begin{pmatrix}
\frac{1}{z_1-y_1} R_{n}(z_1,y_1) & \ldots & \frac{1}{z_1-y_K} R_{n}(z_1,y_K) \\
\vdots & & \vdots \\
\frac{1}{z_L-y_1} R_{n}(z_L,y_1) & \ldots  & \frac{1}{z_L-y_K} R_{n}(z_L,y_K) \\
(Y_{1,1})_{\vecn_{0},\vecm_{0}}(y_1) & \ldots & (Y_{1,1})_{\vecn_{0},\vecm_{0}}(y_K)\\
\vdots & & \vdots \\  (Y_{1,1})_{\vecn_{K-L-1},\vecm_{K-L-1}}(y_1) & \ldots &
(Y_{1,1})_{\vecn_{K-L-1},\vecm_{K-L-1}}(y_K)
\end{pmatrix}.
\end{multline}
Here all $R_{n}$ matrices are taken with respect to the pair of multi-indices
$\vecn_{0},\vecm_{0}$.
\item[(b)] If $L\geq K$ then
\begin{multline}
P_{n}^{[K,L]}(y_1,\ldots,y_K,z_1,\ldots,z_L) = \left(\frac{\prod_{i,j}
(z_j-y_i)}{\prod_{1\leq i<j\leq K} (y_i-y_j)\prod_{1\leq i<j\leq L} (z_j-z_i)}\right)^q\times \\
\det\begin{pmatrix}
\frac{1}{z_1-y_1} L_{n}(y_1,z_1) & \ldots & \frac{1}{z_L-y_1} L_{n}(y_1,z_L) \\
\vdots & & \vdots \\
\frac{1}{z_1-y_K} L_{n}(y_K,z_1) & \ldots  & \frac{1}{z_L-y_K} L_{n}(y_K,z_L) \\
(Y_{2,2})_{\vecn_{0},\vecm_{0}}(z_1) & \ldots & (Y_{2,2})_{\vecn_{0},\vecm_{0}}(z_L)\\
\vdots & & \vdots \\  (Y_{2,2})_{\vecn_{-(L-K-1)},\vecm_{-(L-K-1)}}(z_1) &
\ldots & (Y_{2,2})_{\vecn_{-(L-K-1)},\vecm_{-(L-K-1)}}(z_L)
\end{pmatrix}.
\end{multline}
Here all $L_{n}$ matrices are taken with respect to the pair of multi-indices
$\vecn_{0},\vecm_{0}$, and we assume that $L-K\leq 1+\min_{l=1}^{q} m_l$.
\end{itemize}
\end{theorem}

\subsection{About the proofs of Theorems \ref{theorem:main3} and \ref{theorem:main4}(c)}

Theorems \ref{theorem:main3} and \ref{theorem:main4}(c) will be the key results
from which all other main theorems will follow. The proofs of these two
theorems will be based on a Schur complement formula for the kernel $K_n(x,y)$.
This formula is valid for an arbitrary weight matrix
\begin{equation}\label{weightmatrixgeneral} W(x) =
\begin{pmatrix}W_{k,l}(x)\end{pmatrix}_{k=1,\ldots,p,l=1,\ldots,q}
\end{equation}
of size $p$ by $q$. Thus $W(x)$ does not need to have the rank-one
factorization \eqref{weightmatrixrankone}. For any $k\in\{1,\ldots,p\}$ and
$l\in\{1,\ldots,q\}$ let
\begin{equation*}
\mu_j^{(k,l)} = \int_{-\infty}^{\infty} x^jW_{k,l}(x)\ dx\in\er
\end{equation*}
denote the $j$th moment of the weight function $W_{k,l}(x)$, and for any
$r,s\in\enn$ let
\begin{equation*}
H^{(k,l)}_{r,s} = \begin{pmatrix} \mu^{(k,l)}_{i+j}
\end{pmatrix}_{i=0,\ldots,r-1,j=0,\ldots,s-1}
\end{equation*}
be the $r\times s$ Hankel matrix formed out of these moments. Stack these
matrices in the block Hankel matrix
\begin{equation}\label{blockhankelmatrix}
H_{\vecn,\vecm} = \begin{pmatrix} H^{(1,1)}_{n_1,m_1} &
\ldots & H^{(1,q)}_{n_1,m_q} \\ \vdots & & \vdots \\
H^{(p,1)}_{n_p,m_1} & \ldots & H^{(p,q)}_{n_p,m_q}
\end{pmatrix}.
\end{equation}
Note that this matrix is of size $n\times n$.

Recall that for a block matrix
\begin{equation}\label{defSchurcomp1}
M:=\begin{pmatrix} A & B \\ C & D \end{pmatrix}
\end{equation}
with $A$ square of size $k\times k$ (say), the \emph{Schur complement} of $M$
with respect to the submatrix $D$ is defined as the matrix
\begin{equation}\label{defSchurcomp2} S_{M,D}:=D-CA^{-1}B,\end{equation}
provided $A$ is invertible. Schur complements are also known as
\emph{quasi-determinants} in the literature \cite{Gelfand}. For more
information on Schur complements see Section \ref{subsection:defSchur}.

\begin{proposition}\label{proposition:Schurkernel}
The kernel $K_n(x,y)$ equals the Schur complement of the matrix
\begin{equation}\label{Schurcomplement:kernel} -\left(\begin{array}{ccc|ccc}
H_{n_1,m_1}^{(1,1)} & \ldots & H_{n_1,m_q}^{(1,q)} & \mathbf{y}_{n_1} &  & \mathbf{0}\\
\vdots & & \vdots & & \ddots & \\
H_{n_p,m_1}^{(p,1)} & \ldots & H_{n_p,m_q}^{(p,q)} & \mathbf{0} &  &
\mathbf{y}_{n_p}\\
\hline \vphantom{\overset{O}{\sim}}\mathbf{x}_{m_1}^T &  & \mathbf{0} & 0 & \ldots & 0\\ & \ddots & & \vdots & & \vdots  \\
\mathbf{0} &  & \mathbf{x}_{m_q}^T & 0 & \ldots & 0
\end{array}\right)
\end{equation}
with respect to its bottom right $q\times p$ submatrix. Here the superscript
${}^T$ denotes the matrix transpose, and we use the column vector notations
$\mathbf{y}_{m} = \begin{pmatrix}1 & y & \ldots & y^{m-1}\end{pmatrix}^T$ and
$\mathbf{x}_{m} = \begin{pmatrix}1 & x & \ldots & x^{m-1}\end{pmatrix}^T$.
\end{proposition}

Proposition \ref{proposition:Schurkernel} will be proved in Section
\ref{subsection:kernel}. Variants of this proposition in different contexts can
be found in \cite{Bor,Miranian}, see also \cite{Baik}.

Proposition \ref{proposition:Schurkernel} will allow the quantities $\det
L_n(y,z)$ and $\det R_n(z,y)$ in Theorem \ref{theorem:main3} to be written as a
ratio of two determinants. In Section \ref{subsection:prooftheorem:ratio} we
will expand these determinants using an adaptation of the classical argument of
Heine, and this will lead us to the proof of Theorem~\ref{theorem:main3}. In
Section \ref{subsection:prooftheorem:ratio:more} we will give a similar
argument proving Theorem \ref{theorem:main4}(c).

\subsection{Background on Schur complements}
\label{subsection:defSchur}

Since Schur complements frequently occur in this paper, we find it convenient
to list here some basic properties.

Recall the setting in \eqref{defSchurcomp1} and \eqref{defSchurcomp2}. Schur
complements are related to Gaussian elimination as follows:
\begin{equation}\label{Gaussianelim}
\begin{pmatrix} A & B \\ C & D \end{pmatrix} \sim \begin{pmatrix} A & 0 \\ 0 & S_{M,D}
\end{pmatrix},
\end{equation}
where the $\sim$ symbol relates matrices that can be obtained from one another
by multiplying on the left and on the right with square transformation matrices
of the form $\begin{pmatrix} I_k & 0 \\ * & I
\end{pmatrix}$ and $\begin{pmatrix} I_k & * \\ 0
& I
\end{pmatrix}$, respectively. The procedure to move from the left to the
right hand side of \eqref{Gaussianelim} is sometimes called \emph{Gaussian
elimination with pivot block $A$}, see \cite{GVL}.

In the special case where $D$ and (hence) $S_{M,D}$ are square matrices, one
obtains from \eqref{Gaussianelim} the determinant relation
\begin{equation}\label{Schurdeterminants} \det S_{M,D} = \frac{\det M}{\det
A}.\end{equation}

Let $s_{i,j}$ denote the $(i,j)$th entry of $S = S_{M,D}$. By
\eqref{defSchurcomp2} one has
\begin{equation}\label{Schur:entrywise1}
s_{i,j} = d_{i,j} - \mathbf{c}_{i}A^{-1}\mathbf{b}_{j},
\end{equation}
where $\mathbf{c}_{i}$ denotes the $i$th row of $C$ and $\mathbf{b}_{j}$
denotes the $j$th column of $B$. Now \eqref{Schur:entrywise1} is nothing but
the Schur complement of the matrix
\begin{equation}\label{Schur:entrywise2}
\begin{pmatrix}A & \mathbf{b}_{j}\\
\mathbf{c}_{i} & d_{i,j}\end{pmatrix}\end{equation} with respect to the entry
$d_{i,j}$. From \eqref{Schurdeterminants} we then obtain the determinantal
expression
\begin{equation}\label{Schur:entrywise3}
s_{i,j} = \frac{\det \begin{pmatrix}A & \mathbf{b}_{j}\\
\mathbf{c}_{i} & d_{i,j}\end{pmatrix}}{\det A}.
\end{equation}
This shows that the entries of the Schur complement are ratios of determinants.

More generally than \eqref{Schur:entrywise2}, one may observe that Schur
complements are well-behaved with respect to taking submatrices, in the sense
that any submatrix of the Schur complement \eqref{defSchurcomp2} is itself a
Schur complement, of an appropriate submatrix of \eqref{defSchurcomp1}.

Finally, observe that for any matrix $U$ of appropriate size, the matrix
product $US_{M,D} = U(D-CA^{-1}B)$ can again be written as a Schur complement,
of the matrix
\begin{equation}\label{Schurcomp:multiply}
\begin{pmatrix} A & B \\ UC & UD \end{pmatrix}.
\end{equation}
A similar fact holds for matrix products of the form $S_{M,D}V$.

\subsection{Outline of the rest of the paper}

The rest of this paper is organized as follows. In Section~\ref{section:MOP} we
discuss auxiliary results and notations that are used in the proofs of the main
theorems, and we prove Proposition~\ref{proposition:Schurkernel}.
Theorems~\ref{theorem:main3} and \ref{theorem:main4}(c) are then proved in
Section~\ref{section:proofmaintheorem}, by an adaptation of the classical
argument of Heine. Theorems~\ref{theorem:main1} and \ref{theorem:main2} are
obtained as limiting cases of Theorem~\ref{theorem:main3}. Finally, the
generalizations to arbitrary products and ratios of characteristic polynomials
are proved in Section~\ref{section:proofmaintheorem:ratios}.

\section{Preliminaries for the proofs}
\label{section:MOP}

In this section we collect some auxiliary results and notations that are used
in the proofs of the main theorems. This section is organized as follows.
Vector orthogonal polynomials (which include multiple orthogonal polynomials as
a special case) are defined in Section~\ref{subsection:vectorOP} and their
relation to the RH problem in Section~\ref{subsection:solutionRHP}. The
connection with block Hankel determinants and the duality relations are
discussed in Sections~\ref{subsection:blockHankel} and
\ref{subsection:duality}. In Section~\ref{subsection:kernel} we discuss the
Christoffel-Darboux kernel, leading to the proof of
Proposition~\ref{proposition:Schurkernel}. Finally, the related quantities
$L_n$ and $R_n$ are investigated in Section~\ref{subsection:RnLn}.

\begin{remark}
Practical note: the reader who is mainly interested in Theorems
\ref{theorem:main1}, \ref{theorem:main2} and \ref{theorem:main3} and who wants
to take Proposition \ref{proposition:Schurkernel} for granted may skip this
entire section and move directly to Section~\ref{section:proofmaintheorem}.
\end{remark}

\subsection{Vector orthogonal polynomials}
\label{subsection:vectorOP}

Let $p,q\in\mathbb N $ and consider the multi-indices $\vecn,\vecm$ as in the
beginning of Section \ref{subsection:detpointintro}, but now with $|\vecn| =
|\vecm| +1$. Let $W(x)$ be a $p$ by $q$ weight matrix as in
\eqref{weightmatrixgeneral}. Our definition of vector orthogonal polynomials
will be the same as in Sorokin-Van Iseghem \cite{SVI}; this definition includes
multiple orthogonal polynomials as a special case, when the weight matrix
$W(x)$ has the rank-one factorization \eqref{weightmatrixrankone}.

Let $\Pee_{\vecn}$ be the following space of polynomial vectors
\begin{equation}\label{def:polspace}
\Pee_{\vecn} = \{\vecP(x) = \begin{pmatrix}P_1(x)\\ \vdots \\
P_p(x)\end{pmatrix}\textrm{ with $P_k(x)$ a polynomial of degree at most }
n_k-1 \}.
\end{equation}
The \emph{standard basis} $(\vecA_1(x),\ldots,\vecA_{|\vecn|}(x))$ of the
vector space $\Pee_{\vecn}$ is defined by the column vectors of the $p$ by
$|\vecn|$ matrix
\begin{equation}\label{standardbasis:A}
\begin{pmatrix} \vecA_1(x) & \ldots & \vecA_{|\vecn|}(x) \end{pmatrix}
= \left(\begin{array}{cccc|cccc|cc} 1 & x & \ldots & x^{n_1-1} & 0 & 0 & \ldots & 0 & \ldots & 0 \\
0 & 0 & \ldots & 0 & 1 & x & \ldots & x^{n_2-1} & \ldots & 0 \\
\vdots & \vdots & & \vdots & \vdots & \vdots & & \vdots & & \vdots \\
0 & 0 & \ldots & 0 & 0 & 0 & \ldots & 0 & \ldots & x^{n_p-1}
 \end{array}\right).
\end{equation}
Similarly, let $\Pee_{\vecm}$ denote the vector space
\begin{equation}\label{def:polspace:Q}
\Pee_{\vecm} = \{\vecQ(x) = \begin{pmatrix}Q_1(x)\\ \vdots \\
Q_q(x)\end{pmatrix}\textrm{ with $Q_l(x)$ a polynomial of degree at most }
m_l-1 \},
\end{equation}
and define the standard basis $(\vecB_1(x),\ldots,\vecB_{|\vecm|}(x))$ for
$\Pee_{\vecm}$ by the columns of the $q$ by $|\vecm|$ matrix
\begin{equation}\label{standardbasis:B}
\begin{pmatrix} \vecB_1(x) & \ldots & \vecB_{|\vecm|}(x) \end{pmatrix}
= \left(\begin{array}{cccc|cccc|cc} 1 & x & \ldots & x^{m_1-1} & 0 & 0 & \ldots & 0 & \ldots & 0 \\
0 & 0 & \ldots & 0 & 1 & x & \ldots & x^{m_2-1} & \ldots & 0 \\
\vdots & \vdots & & \vdots & \vdots & \vdots & & \vdots & & \vdots \\
0 & 0 & \ldots & 0 & 0 & 0 & \ldots & 0 & \ldots & x^{m_q-1}
 \end{array}\right).
\end{equation}
Note that we use boldface notation to denote vector-valued objects.

We say that $\vecP(x) = \vecP_{\vecn,\vecm}(x)\in\Pee_{\vecn}$ is a
\emph{vector orthogonal polynomial} with respect to the multi-indices $\vecn$,
$\vecm$ and the weight matrix $W(x)$ if
$$ \int_{-\infty}^{\infty} \vecP_{\vecn,\vecm}^T(x)W(x)\vecQ(x) dx = 0$$
for every $\vecQ\in\Pee_{\vecm}$. To stress the dependence on the weight matrix
we will sometimes write $\vecP_{\vecn,\vecm}(x) = \vecP_{\vecn,\vecm,W}(x)$.

The coefficients of the vector orthogonal polynomials can be found from a
homogeneous linear system with $|\vecn|$ unknowns (polynomial coefficients) and
$|\vecm|$ equations (orthogonality conditions).
The restriction $|\vecn|=|\vecm|+1$ guarantees that this system has a
nontrivial solution. If $\vecP_{\vecn,\vecm}(x)$ is unique up to a
multiplicative factor then the pair of multi-indices $\vecn,\vecm$ is called
\emph{normal}.

Assume that $\vecn,\vecm$ is a normal pair of indices. Then
$\vecP_{\vecn,\vecm}(x)$ is said to satisfy \begin{itemize}\item the
\emph{normalization of type I with respect to the $l$th index},
$l\in\{1,\ldots,q\}$, if
\begin{equation*}\label{type1norm}
\int_{-\infty}^{\infty} \vecP^T_{\vecn,\vecm}(x) W(x)\vecQ(x)\ dx = 1
\end{equation*}
where $\vecQ(x)$ is the vector consisting of zeros except for the $l$th entry
which is $x^{m_l}$.
\item the \emph{normalization of type II with respect to the $k$th
index}, $k\in\{1,\ldots,p\}$, if the $k$th component $(P_k)_{\vecn,\vecm}(x)$
is monic, i.e., if
\begin{equation*}\label{type2norm}
(P_k)_{\vecn,\vecm}(x) = x^{n_k-1}+O(x^{n_k-2}).
\end{equation*}
\end{itemize}

The vector orthogonal polynomials $\vecP_{\vecn,\vecm}(x)$ corresponding to the
above normalizations, if they exist, will be denoted as
$\vecP_{\vecn,\vecm}^{(I,l)}(x)$ and $\vecP_{\vecn,\vecm}^{(II,k)}(x)$,
respectively.


\subsection{Solution to the RH problem}
\label{subsection:solutionRHP}

Consider again the RH problem in Section \ref{subsection:RHintro} with $|\vecn|
= |\vecm|$ but now with $W(x)$ in \eqref{weightmatrixrankone} replaced by an
arbitrary matrix \eqref{weightmatrixgeneral}. The solution to the RH problem,
if it exists, is uniquely described by vector orthogonal polynomials and their
Cauchy transforms \cite{DK2}. More precisely, consider again the partition of
$Y(z)$ as in \eqref{partitionY}. Then one has that
\begin{itemize}
\item $Y_{1,1}(z)$ has its $k$th row given by the row vector
\begin{equation}\label{MOPtype2}
(\vecP^{(II,k)}_{\vecn+\vece_k,\vecm})^T(z),
\end{equation}
$k=1,\ldots,p$.
\item $Y_{2,1}(z)$ has its $l$th row given by
\begin{equation}\label{MOPtype1}
-2\pi i(\vecP^{(I,l)}_{\vecn,\vecm-\vece_l})^T(z),
\end{equation}
$l=1,\ldots,q$.
\item $Y_{1,2}(z)$ has its $k$th row given by
\begin{equation}\label{cauchytransform1}
-\frac{1}{2\pi i}\int_{-\infty}^{\infty}
\frac{(\vecP^{(II,k)}_{\vecn+\vece_k,\vecm})^T(x)W(x)}{z-x}\ dx,
\end{equation}
$k=1,\ldots,p$.
\item $Y_{2,2}(z)$ has its $l$th row given by
\begin{equation}\label{cauchytransform2}
\int_{-\infty}^{\infty}
\frac{(\vecP^{(I,l)}_{\vecn,\vecm-\vece_l})^T(x)W(x)}{z-x}\ dx,\end{equation}
$l=1,\ldots,q$.
\end{itemize}
Here we are using the notation $\vece_l$ to denote the standard basis vector
which is zero except for its $l$th entry which is one. The length of this
vector will always be clear from the context.

For example, in the special case where $p=q=2$, the solution $Y(z)$ is given by
the $4\times 4$ matrix
\begin{equation*}\label{RHmatrix4x40} Y(z) =
D\times\left(\begin{array}{cccc}
(P_1)_{\vecn+\vece_1,\vecm}^{(II,1)} & (P_2)_{\vecn+\vece_1,\vecm}^{(II,1)} & * & * \\
(P_1)_{\vecn+\vece_2,\vecm}^{(II,2)} & (P_2)_{\vecn+\vece_2,\vecm}^{(II,2)} & * & * \\
(P_1)_{\vecn,\vecm-\vece_1}^{(I,1)} & (P_2)_{\vecn,\vecm-\vece_1}^{(I,1)} & * & * \\
(P_1)_{\vecn,\vecm-\vece_2}^{(I,2)} & (P_2)_{\vecn,\vecm-\vece_2}^{(I,2)} & * & * \\
\end{array}\right),
\end{equation*}
where $D := \diag(1,1,-2\pi i,-2\pi i)$, and where the entries denoted with $*$
are Cauchy transforms as in \eqref{cauchytransform1} and
\eqref{cauchytransform2}.

\subsection{Biorthogonality and moment matrix}
\label{subsection:blockHankel}

Recall the standard bases $\vec{A}_i(x)$ and $\vec{B}_j(x)$ in
\eqref{standardbasis:A} and \eqref{standardbasis:B}, where now $|\vecn| =
|\vecm| =:n$. The \emph{moment matrix} $M$ with respect to these bases is of
size $n$ by $n$ and has entries
\begin{equation*}
M_{i,j} = \int_{-\infty}^{\infty} \vecA_i^T(x) W(x)\vecB_j(x)\ dx.
\end{equation*}
The moment matrix coincides with the block Hankel matrix $H_{\vecn,\vecm}$ in
\eqref{blockhankelmatrix}.

Consider ordered bases $(\vecP_1(x),\ldots,\vecP_n(x))$ and
$(\vecQ_1(x),\ldots,\vecQ_n(x))$ of the vector spaces $\Pee_{\vecn}$ and
$\Pee_{\vecm}$, respectively. We say that these bases are \emph{biorthogonal}
if
\begin{equation}\label{biorthogonality}
\int_{-\infty}^{\infty} \vecP_i^T(x) W(x)\vecQ_j(x)\ dx = \delta_{i,j},
\end{equation}
the Kronecker delta.

\begin{lemma}\label{RHP:solvable:Hankel}
Assume that $|\vecn|=|\vecm|$. Then the following statements are equivalent
\begin{enumerate}
\item The RH problem for $Y(z)$ in Section \ref{subsection:RHintro} is solvable.
\item There is no non-zero vector in $\Pee_{\vecn}$ which is biorthogonal to
the entire space $\Pee_{\vecm}$.
\item $\det H_{\vecn,\vecm}\neq 0$.
\item There exist biorthogonal bases for the spaces $\Pee_{\vecn}$ and
$\Pee_{\vecm}$.
\end{enumerate}
\end{lemma}

\bewijs. First consider the equivalence between Statements 1 and 2. Suppose
that statement~2 holds true. Then each of the polynomial vectors in
\eqref{MOPtype2}--\eqref{MOPtype1} exists \cite{DK2} and therefore the RH
problem for $Y(z)$ is solvable. Conversely, suppose that statement~2 does not
hold. Then the polynomial vectors in \eqref{MOPtype2}--\eqref{MOPtype1}, if
they exist, are not unique and therefore the RH problem for $Y(z)$ cannot be
solvable, since that would contradict the uniqueness of $Y(z)$. Finally, the
equivalence between Statements 2--4 follows from standard linear algebra
arguments whose description we omit. $\bol$


\subsection{Duality}
\label{subsection:duality}

The role of the biorthogonal bases $\vecP_i$ and $\vecQ_j$ is interchanged by
swapping the multi-indices $\vecn$ and $\vecm$ and by transposing the weight
matrix $W(x)$:
\begin{equation}\label{duality} (\vecn,\vecm,W(x))\leftrightarrow (\vecm,\vecn,W^T(x)).\end{equation}
This duality has a nice form in terms of the RH problem \cite{AvMV,DK2}. If we
partition the RH matrices corresponding to the original and dual data in
\eqref{duality} as
\begin{equation*} Y_{\vecn,\vecm,W} =
\begin{array}{ll} & \begin{array}{ll} \ \ p & \ \ \ \ q \end{array}\\
\begin{array}{l} p \\ q\end{array} \hspace{-4mm} & \begin{pmatrix}Y_{1,1} & Y_{1,2}\\ Y_{2,1} &
Y_{2,2}
\end{pmatrix},\end{array}\qquad Y_{\vecm,\vecn,W^T} =
\begin{array}{ll} & \begin{array}{ll} \ \ q & \ \ \ \ p \end{array}\\
\begin{array}{l} q \\ p\end{array} \hspace{-4mm} & \begin{pmatrix}\tilde Y_{1,1} & \tilde Y_{1,2}\\ \tilde Y_{2,1} &
\tilde Y_{2,2}
\end{pmatrix},\end{array}
\end{equation*}
respectively, then it holds that
\begin{equation}\label{dual:partitioning}
\begin{pmatrix} \tilde Y_{1,1} & \tilde Y_{1,2} \\ \tilde Y_{2,1} & \tilde Y_{2,2} \end{pmatrix}
= \begin{pmatrix} Y_{2,2} & -Y_{2,1} \\ -Y_{1,2} & Y_{1,1} \end{pmatrix}^{-T},
\end{equation}
where the superscript ${}^{-T}$ denotes the inverse transpose. From a
well-known theorem on minors of the inverse matrix, see e.g.\ \cite[page.~21,
eq.~(33)]{Gant}, and the general fact that $\det Y(z)\equiv 1$ it then follows
that
\begin{equation}\label{jacobi1} \det \tilde Y_{1,1}(y) = \det Y_{1,1}(y)
 \end{equation} and
 \begin{equation}\label{jacobi2}
\det \tilde Y_{2,2}(z) = \det Y_{2,2}(z).
\end{equation}
The duality relations \eqref{jacobi1}--\eqref{jacobi2} take care of a point
made earlier, see the remark after the statement of Theorem
\ref{theorem:main2}.

From another application of \cite[loc. cit.]{Gant} it also follows that
\begin{equation}\label{jacobi3}
\det\left(\begin{pmatrix}I_p & 0 \end{pmatrix}Y^{-1}(z)Y(y)\begin{pmatrix}I_p \\
0
\end{pmatrix}\right) = \det\left(\begin{pmatrix}0 & I_q
\end{pmatrix}Y^{-1}(y)Y(z)\begin{pmatrix}0 \\ I_q
\end{pmatrix}\right).
\end{equation}
Comparing this with \eqref{def:Ln} and \eqref{def:Rn}, this establishes the
fact that $\det R_n(z,y) = \det L_n(y,z)$, which is of course implicit in the
statement of Theorem \ref{theorem:main3}.

\subsection{The Christoffel-Darboux kernel}
\label{subsection:kernel}

In this section we prove the Schur complement formula for the kernel $K_n(x,y)$
in Proposition~\ref{proposition:Schurkernel}. To this end we identify
$K_n(x,y)$ as a reproducing kernel. In this way we will show that
\eqref{def:kernel:block} corresponds to the kernels familiar from the theory of
matrix orthogonal polynomials \cite{Miranian} and vector orthogonal
polynomials~\cite{SVI}.

We will assume throughout that the RH problem for $Y(z)$ is solvable, which is
of course necessary in order for \eqref{def:kernel:block} to make sense.

\begin{lemma}
The kernel $K_n(x,y)$ is a $q\times p$ matrix whose $(i,j)$ entry is a
bivariate polynomial of degree at most $m_i-1$ in $x$ and $n_j-1$ in $y$.
Equivalently,
\begin{equation}\label{kernel:bivariate}
K_n(x,y) = \sum_{i,j=1}^{n} c_{i,j}\vecB_i(x)\vecA^T_j(y),
\end{equation}
for certain $c_{i,j}$, where $\vecA_j$ and $\vecB_i$ denote the basis vectors
in \eqref{standardbasis:A} and \eqref{standardbasis:B}.
\end{lemma}

\bewijs. The expression $\begin{pmatrix}0 & I_q \end{pmatrix} Y^{-1}(x)Y(y)
\begin{pmatrix}I_p \\ 0 \end{pmatrix}$ in the numerator of \eqref{def:kernel:block}
is made out of the
polynomial entries of the matrices $Y^{-1}(x)$ and $Y(y)$, cf.\
\eqref{MOPtype2}--\eqref{MOPtype1} and \eqref{dual:partitioning}. Since this
expression vanishes if $x=y$, each of these polynomials is divisible by $y-x$,
so the denominator of \eqref{def:kernel:block} can be divided out.$\bol$

It was shown in \cite{DK2} that for a rank-one weight matrix $W(x)$ as in
\eqref{weightmatrixrankone}, the scalar-valued kernel $\hat K_n(x,y)$ in
\eqref{def:kernel:scalar} has a certain reproducing property. We now obtain a
similar property for the matrix-valued kernel $K_n(x,y)$ (without restrictions
on $W(x)$).

\begin{proposition}\label{proposition:reproducing}
The kernel $K_n$ satisfies the reproducing property
\begin{equation}\label{reproducing}
\int_{-\infty}^{\infty} K_n(x,y) W(y)\vecQ(y)\ dy = \vecQ(x)
\end{equation}
for any $x\in\cee$ and $\vecQ\in\Pee_{\vecm}$. Moreover, this reproducing
property uniquely characterizes $K_n(x,y)$ in the class of bivariate polynomial
matrices of the form \eqref{kernel:bivariate}.
\end{proposition}

\bewijs. First we establish \eqref{reproducing}. This can be done by adapting
the proof in \cite{DK2}. However, we give a more streamlined proof. We have
$$ \int K_n(x,y) W(y)\vecQ(y)\ dy =  \int K_n(x,y) W(y)(\vecQ(y)-\vecQ(x))\ dy
+ \int K_n(x,y) W(y)\vecQ(x)\ dy. $$ By inserting \eqref{def:kernel:block} this
becomes
\begin{equation}\label{proof:reproducing:twoterms} = \begin{pmatrix}0 & I_q
\end{pmatrix} Y^{-1}(x)
\begin{pmatrix}\frac{1}{2\pi
i}\int Y_{1,1}(y)W(y)\frac{\vecQ(y)-\vecQ(x)}{y-x}\ dy \\ \frac{1}{2\pi i}\int
Y_{2,1}(y)W(y)\frac{\vecQ(y)-\vecQ(x)}{y-x}\ dy \end{pmatrix}+
\begin{pmatrix}0 & I_q \end{pmatrix} Y^{-1}(x)
\begin{pmatrix}\frac{1}{2\pi i}\int \frac{Y_{1,1}(y)W(y)}{y-x}\ dy \\ \frac{1}{2\pi i}\int \frac{Y_{2,1}(y)W(y)}{y-x}\ dy \end{pmatrix}\vecQ(x).
\end{equation}
The first term in \eqref{proof:reproducing:twoterms} is zero since
$\frac{\vecQ(y)-\vecQ(x)}{y-x}$ is a polynomial vector in $y$ whose $l$th entry
has degree at most $m_l-2$, $l=1,\ldots,q$, and by invoking the orthogonality
relations of the vector orthogonal polynomials. For the second term in
\eqref{proof:reproducing:twoterms}, we recognize the defining relations for the
Cauchy transforms in the second block column of the RH problem, see
\eqref{cauchytransform1}--\eqref{cauchytransform2}, so we get
$$ = \begin{pmatrix}0 & I_q \end{pmatrix} Y^{-1}(x)
\begin{pmatrix}Y_{1,2}(x) \\ Y_{2,2}(x) \end{pmatrix}\vecQ(x) =
\begin{pmatrix}0 & I_q \end{pmatrix}\begin{pmatrix}0 \\ I_q \end{pmatrix}\vecQ(x) = \vecQ(x).
$$
This establishes \eqref{reproducing} when $x\in\cee\setminus\er$; the case
where $x\in\er$ follows by continuity.

Next we show that the reproducing property \eqref{reproducing} uniquely
characterizes $K_n(x,y)$ in the class of bivariate polynomial matrices
\eqref{kernel:bivariate}. By Lemma \ref{RHP:solvable:Hankel}, $1\Rightarrow 4$,
we can choose biorthogonal bases for the polynomial vector spaces
$\Pee_{\vecn}$ and $\Pee_{\vecm}$; let $(\vecP_k)_{k=1}^n$ and
$(\vecQ_k)_{k=1}^n$ be such bases. Any $q$ by $p$ matrix $M_n(x,y)$ of the form
\eqref{kernel:bivariate} can be rewritten as
$$ M_n(x,y) = \sum_{i,j=1}^{n} \tilde c_{i,j}\vecQ_i(x)\vecP^T_j(y)
$$
for suitable constants $\tilde c_{i,j}$. But then the reproducing property
\eqref{reproducing} and the biorthogonality relations \eqref{biorthogonality}
imply that $\tilde c_{i,j}=\delta_{i,j}$, the Kronecker delta. This ends the
proof. $\bol$

Proposition \ref{proposition:reproducing} has the following dual version.

\begin{proposition}\label{proposition:reproducing:dual}
The kernel $K_n$ satisfies the reproducing property
\begin{equation}\label{reproducing:dual}
\int_{-\infty}^{\infty} \vecP^T(x)W(x)K_n(x,y)\ dx = \vecP^T(y)
\end{equation}
for any $y\in\cee$ and $\vecP\in\Pee_{\vecn}$. Moreover, this reproducing
property uniquely characterizes $K_n(x,y)$ in the class of bivariate polynomial
matrices of the form \eqref{kernel:bivariate}.
\end{proposition}

Now we are ready for the

\textsc{Proof of Proposition \ref{proposition:Schurkernel}}. Denote with
$M_n(x,y)$ the Schur complement of \eqref{Schurcomplement:kernel}. We check the
reproducing kernel property
\begin{equation}\label{reproducing1} \int \vecA^T(x)W(x)M_n(x,y)\ dx = \vecA^T(y), \end{equation}
for each column vector of the form
$\vecA(x) = \begin{pmatrix} 0 & \ldots & x^{i_k} & \ldots & 0
\end{pmatrix}^T$ with $i_k\in\{0,\ldots,n_k-1\}$, $k\in\{1,\ldots,p\}$.
(These vectors $\vecA(x)$ are the standard basis of the space $\Pee_{\vecn}$.)
By using the multiplication property \eqref{Schurcomp:multiply} for Schur
complements, the left hand side of \eqref{reproducing1} equals the Schur
complement of
\begin{equation}\label{reproducing:proof1}
-\left(\begin{array}{ccc|ccc}
H_{n_1,m_1}^{(1,1)} & \ldots & H_{n_1,m_q}^{(1,q)} & \mathbf{y}_{n_1} &  & \mathbf{0}\\
\vdots & & \vdots & & \ddots & \\
H_{n_p,m_1}^{(p,1)} & \ldots & H_{n_p,m_q}^{(p,q)} & \mathbf{0} &  &
\mathbf{y}_{n_p}\\
\hline \vphantom{\overset{O}{\sim}}\mathbf{h}_{i_k,m_1}^{(k,1)} &  &
\mathbf{h}_{i_k,m_q}^{(k,q)} & 0 & \ldots & 0
\end{array}\right)
\end{equation}
with respect to its bottom right $1\times p$ submatrix, where we use the row
vector notation $\mathbf{h}_{i_k,m_l}^{(k,l)} = \begin{pmatrix}
\mu^{(k,l)}_{i_k+j}
\end{pmatrix}_{j=0,\ldots,m_l-1}$. Subtracting
from the last row of \eqref{reproducing:proof1} the $i_k$th row from the $k$th
block row, we are led to the Schur complement of
\begin{equation}\label{reproducing:proof2}
-\left(\begin{array}{ccc|ccccc}
H_{n_1,m_1}^{(1,1)} & \ldots & H_{n_1,m_q}^{(1,q)} & \mathbf{y}_{n_1} & & & & \mathbf{0}\\
\vdots & & \vdots & & & \ddots & & \\
H_{n_p,m_1}^{(p,1)} & \ldots & H_{n_p,m_q}^{(p,q)} & \mathbf{0} & & & &
\mathbf{y}_{n_p}\\
\hline \vphantom{\overset{O}{\sim}} \mathbf{0} &  & \mathbf{0} & 0 & \ldots &
-y^{i_k} & \ldots & 0
\end{array}\right),
\end{equation}
where now the last row is zero except for the entry $-y^{i_k}$ in the $k$th
column of the second block column. But clearly, the Schur complement of
\eqref{reproducing:proof2} is just the row vector $\vecA^T(y) =
\begin{pmatrix} 0 & \ldots & y^{i_k} & \ldots & 0
\end{pmatrix}$. This establishes the reproducing kernel property \eqref{reproducing1}.

To finish the proof of the proposition, we note that the degree structure of
\eqref{Schurcomplement:kernel} implies that $M_n(x,y)$ has the form in
\eqref{kernel:bivariate}. Thus the proof is ended by invoking the uniqueness
part of Proposition \ref{proposition:reproducing:dual}.  $\bol$

\subsection{The matrices $L_n$ and $R_n$}
\label{subsection:RnLn}

In this section we consider in more detail the quantities $L_n$ and $R_n$ which
are derived from the kernel $K_n$. Let us first establish the equivalence
between the two different formulae in the definition of $L_n$ in
\eqref{def:Ln}. From the identity $y-x = (y-z)+(z-x)$ we obtain
\begin{eqnarray}
\nonumber \int \frac{y-x}{z-x}K_n(y,x)W(x)\ dx &=& \int K_n(y,x)W(x)I_q\ dx+\int \frac{y-z}{z-x}K_n(y,x)W(x)\ dx\\
\label{proof:vanishing0} &=& I_q-(z-y)\int K_n(y,x)W(x)/(z-x)\ dx
\end{eqnarray}
by virtue of the reproducing property in Proposition
\ref{proposition:reproducing}. On the other hand, the equality
\begin{equation*}
\int \frac{y-x}{z-x}K_n(y,x)W(x)\ dx = \begin{pmatrix}0 & I_q \end{pmatrix}Y^{-1}(y)Y(z)\begin{pmatrix}0 \\
I_q
\end{pmatrix}
\end{equation*}
follows from \eqref{def:kernel:block} and
\eqref{cauchytransform1}--\eqref{cauchytransform2}. This establishes the
required equality in \eqref{def:Ln}.

The equivalence between the different formulae in \eqref{def:Rn} can be
established similarly.

In section \ref{section:proofmaintheorem:ratios}, we will need the following
property of $L_n$.

\begin{proposition}\label{proposition:vanishingprop}
The matrix $L_n$ in \eqref{def:Ln} satisfies the \lq vanishing property\rq\
\begin{equation}\label{vanishing:dual}
\int_{-\infty}^{\infty} \vecP^T(y)W(y)\frac{L_n(y,z)}{z-y} \ dy = 0
\end{equation}
for any $z\in\cee\setminus\er$ and $\vecP\in\Pee_{\vecn}$.
\end{proposition}

\bewijs. From the expression of $L_n$ in \eqref{proof:vanishing0}, we obtain
\begin{equation*}\label{proof:vanishing1}
\int \vecP^T(y)W(y)\frac{L_n(y,z)}{z-y}\ dy = \int \vecP^T(y)\frac{W(y)}{z-y} \
dy - \int\!\int \frac{\vecP^T(y)W(y)K_n(y,x)W(x)}{z-x}\ dx\ dy
\end{equation*}
which is
$$ = \int \vecP^T(y)\frac{W(y)}{z-y} \
dy - \int \left(\int \vecP^T(y)W(y)K_n(y,x) dy\right)\frac{W(x)}{z-x} dx.
$$
By inserting the reproducing property this becomes
$$ = \int \vecP^T(y)\frac{W(y)}{z-y} \
dy - \int \vecP^T(x)\frac{W(x)}{z-x} \ dx = 0.
$$
$\bol$

Here is the dual version of Proposition \ref{proposition:vanishingprop}.

\begin{proposition}\label{proposition:vanishingprop:dual}
The matrix $R_n$ in \eqref{def:Rn} satisfies the \lq vanishing property\rq\
\begin{equation}\label{vanishing}
\int_{-\infty}^{\infty} \frac{R(z,y)}{z-y} W(y)\vecQ(y)\ dy = 0
\end{equation}
for any $z\in\cee\setminus\er$ and $\vecQ\in\Pee_{\vecm}$.
\end{proposition}

\section{Proofs, part 1}
\label{section:proofmaintheorem}

In this section we prove Theorems \ref{theorem:main1}, \ref{theorem:main2},
\ref{theorem:main3} and \ref{theorem:main4}(c).

\subsection{Proof of Theorem \ref{theorem:main3}}
\label{subsection:prooftheorem:ratio}

First we prove Theorem \ref{theorem:main3}. The proof of this theorem will
follow by expanding the moment determinant in Proposition
\ref{proposition:Schurkernel} in a similar way as in the classical argument of
Heine.

We will restrict ourselves to the proof of \eqref{mainR3a}; the proof of
\eqref{mainR3b} can then be devised in a similar way, or by simply invoking
\eqref{jacobi3}. For notational convenience and to keep things readable, we
give the proof for the case $p=q=2$. The case of general $p,q$ is completely
similar except that it requires more notational burden.

First we expand the left hand side of \eqref{mainR3a}. By definition, the
average ratio of characteristic polynomials $P_n^{[1,1]}(y,z)$ is given by
\begin{multline}\label{proof3:startingchar}
P_n^{[1,1]}(y,z) = \frac{1}{Z_n}\int\ldots\int\prod_{j=1}^n \frac{y-x_j}{z-x_j}\\
\left|\begin{array}{ccc} w_{1,1}(x_1) & \ldots & w_{1,1}(x_n) \\
\vdots & & \vdots \\
x_1^{n_1-1}w_{1,1}(x_1) & \ldots & x_n^{n_1-1}w_{1,1}(x_n) \\
\hline w_{1,2}(x_1) & \ldots & w_{1,2}(x_n) \\
\vdots & & \vdots \\
x_1^{n_2-1}w_{1,2}(x_1) & \ldots & x_n^{n_2-1}w_{1,2}(x_n)
\end{array}\right|\times
\left|\begin{array}{ccc} w_{2,1}(x_1) & \ldots & w_{2,1}(x_n) \\
\vdots & & \vdots \\
x_1^{m_1-1}w_{2,1}(x_1) & \ldots & x_n^{m_1-1}w_{2,1}(x_n) \\
\hline w_{2,2}(x_1) & \ldots & w_{2,2}(x_n) \\
\vdots & & \vdots \\
x_1^{m_2-1}w_{2,2}(x_1) & \ldots & x_n^{m_2-1}w_{2,2}(x_n)
\end{array}\right|\prod_{j=1}^n dx_j.
\end{multline}
We can expand the second determinant according to the Lagrange expansion with
respect to the first block row:
\begin{multline}\label{proof3:Lagrange}
\left|\begin{array}{ccc} w_{2,1}(x_1) & \ldots & w_{2,1}(x_n) \\
\vdots & & \vdots \\
x_1^{m_1-1}w_{2,1}(x_1) & \ldots & x_n^{m_1-1}w_{2,1}(x_n) \\
\hline w_{2,2}(x_1) & \ldots & w_{2,2}(x_n) \\
\vdots & & \vdots \\
x_1^{m_2-1}w_{2,2}(x_1) & \ldots & x_n^{m_2-1}w_{2,2}(x_n)
\end{array}\right| \\
= \sum_{S} (-1)^{S}\begin{vmatrix} w_{2,1}(x_{s_1}) & \ldots & w_{2,1}(x_{s_{m_1}}) \\
\vdots & & \vdots \\
x_{s_1}^{m_1-1}w_{2,1}(x_{s_1}) & \ldots &
x_{s_{m_1}}^{m_1-1}w_{2,1}(x_{s_{m_1}})\end{vmatrix}\begin{vmatrix} w_{2,2}(x_{\bar{s}_1}) & \ldots & w_{2,2}(x_{\bar{s}_{m_2}}) \\
\vdots & & \vdots \\
x_{\bar s_1}^{m_2-1}w_{2,2}(x_{\bar s_1}) & \ldots & x_{\bar
s_{m_2}}^{m_2-1}w_{2,2}(x_{\bar s_{m_2}})\end{vmatrix}
\end{multline}
where the sum is over all $\binom{n}{m_1}$ subsets $S\subset\{1,\ldots,n\}$
with $|S| = m_1$, where we write the elements of $S$ in increasing order:
$s_1<\ldots<s_{m_1}$, similarly for those of the complement $\bar{S} =
\{1,\ldots,n\}\setminus S$, $\bar{s}_1<\ldots<\bar{s}_{m_2}$, and where we
denote by $(-1)^{S}$ the sign of the permutation $s_1,\ldots,s_{m_1},\bar
s_1,\ldots,\bar s_{m_2}$.

For each term of \eqref{proof3:Lagrange} we can incorporate the factor
$(-1)^{S}$ into the other determinant of \eqref{proof3:startingchar} by
permuting its columns. Next we can relabel the variables $x_1,\ldots,x_n$ back
in their original form, and so \eqref{proof3:startingchar} reduces to
\begin{multline*}
P_n^{[1,1]}(y,z) = \frac{\binom{n}{m_1}}{Z_n}\int\ldots\int\prod_{j=1}^n
\frac{y-x_j}{z-x_j}
\left|\begin{array}{ccc} w_{1,1}(x_1) & \ldots & w_{1,1}(x_n) \\
\vdots & & \vdots \\
x_1^{n_1-1}w_{1,1}(x_1) & \ldots & x_n^{n_1-1}w_{1,1}(x_n) \\
\hline w_{1,2}(x_1) & \ldots & w_{1,2}(x_n) \\
\vdots & & \vdots \\
x_1^{n_2-1}w_{1,2}(x_1) & \ldots & x_n^{n_2-1}w_{1,2}(x_n)
\end{array}\right|\\ \begin{vmatrix} w_{2,1}(x_{1}) & \ldots & w_{2,1}(x_{m_1}) \\
\vdots & & \vdots \\
x_1^{m_1-1}w_{2,1}(x_{1}) & \ldots &
x_{m_1}^{m_1-1}w_{2,1}(x_{m_1})\end{vmatrix}
\begin{vmatrix} w_{2,2}(x_{m_1+1}) & \ldots & w_{2,2}(x_{n}) \\
\vdots & & \vdots \\
x_{m_1+1}^{m_2-1}w_{2,2}(x_{m_1+1}) & \ldots &
x_n^{m_2-1}w_{2,2}(x_{n})\end{vmatrix}\prod_{j=1}^n dx_j
\end{multline*}
which on evaluating the Vandermonde determinants becomes
\begin{multline}\label{proof3:lefthandside:reduced}
P_n^{[1,1]}(y,z) = \frac{\binom{n}{m_1}}{Z_n}\int\ldots\int\ \prod_{j=1}^{n}
\frac{y-x_j}{z-x_j}\ \left|\begin{array}{ccc} w_{1,1}(x_1) & \ldots &
w_{1,1}(x_{n}) \\
\vdots  & & \vdots  \\
x_1^{n_1-1} w_{1,1}(x_1)  & \ldots  &
x_{n}^{n_1-1}w_{1,1}(x_{n})  \\
\hline w_{1,2}(x_1)  & \ldots & w_{1,2}(x_{n}) \\
\vdots  & & \vdots  \\
x_1^{n_2-1} w_{1,2}(x_1)  & \ldots & x_{n}^{n_2-1}w_{1,2}(x_{n})
\end{array}\right|\\ \prod_{1\leq i<j\leq m_1} (x_j-x_i)
\prod_{m_1+1\leq i<j\leq n} (x_j-x_i) \prod_{j=1}^{m_1}w_{2,1}(x_j)
\prod_{j=m_1+1}^{n} w_{2,2}(x_j)\prod_{j=1}^n dx_j.
\end{multline}

Next we expand the right hand side of \eqref{mainR3a}. We recall from
Proposition \ref{proposition:Schurkernel} that $-K_n(x,y)$ is the Schur
complement of the matrix
\begin{equation}\label{proof3:startingmatrix0}
\left(\begin{array}{cc|cc} H_{n_1,m_1}^{(1,1)} & H_{n_1,m_2}^{(1,2)} &
\mathbf{y}_{n_1} & \mathbf{0}  \\
H_{n_2,m_1}^{(2,1)} & H_{n_2,m_2}^{(2,2)} &
\mathbf{0} & \mathbf{y}_{n_2} \\
\hline \vphantom{\overset{O}{\sim}}\mathbf{x}^T_{m_1} & 0 & 0 & 0 \\
0 & \mathbf{x}^T_{m_2} & 0 & 0
\end{array}\right).
\end{equation}
This implies that $\frac{1}{z-y}I_2-\int \frac{W(x)K_n(x,y)}{z-x}\ dx$ is the
Schur complement of the matrix
\begin{equation}\label{proof3:startingmatrix0bis}
\left(\begin{array}{cc|cc} H_{n_1,m_1}^{(1,1)} & H_{n_1,m_2}^{(1,2)} &
\mathbf{y}_{n_1} & \mathbf{0}  \\
H_{n_2,m_1}^{(2,1)} & H_{n_2,m_2}^{(2,2)} &
\mathbf{0} & \mathbf{y}_{n_2} \\
\hline \vphantom{\overset{O}{\sim}}\int\mathbf{x}^T_{m_1}W_{1,1}(x)/(z-x)\ dx & \int\mathbf{x}^T_{m_2}W_{1,2}(x)/(z-x)\ dx & 1/(z-y) & 0 \\
\int\mathbf{x}^T_{m_1}W_{2,1}(x)/(z-x)\ dx &
\int\mathbf{x}^T_{m_2}W_{2,2}(x)/(z-x)\ dx & 0 & 1/(z-y)
\end{array}\right).
\end{equation}
From \eqref{Schurdeterminants} this implies the determinant representation
\begin{multline}\label{proof3:startingmatrix1}
\det\left(\frac{1}{z-y}I_2-\int \frac{W(x)K_n(x,y)}{z-x}\ dx\right) \\=
\frac{1}{\kappa}\begin{vmatrix} H_{n_1,m_1}^{(1,1)} & H_{n_1,m_2}^{(1,2)} &
\mathbf{y}_{n_1} & \mathbf{0}  \\
H_{n_2,m_1}^{(2,1)} & H_{n_2,m_2}^{(2,2)} &
\mathbf{0} & \mathbf{y}_{n_2} \\
\int\mathbf{x}_{m_1}^T W_{1,1}(x)/(z-x)\ dx & \int\mathbf{x}_{m_2}^T W_{1,2}(x)/(z-x)\ dx & 1/(z-y) & 0 \\
\int\mathbf{x}_{m_1}^T W_{2,1}(x)/(z-x)\ dx & \int\mathbf{x}_{m_2}^T
W_{2,2}(x)/(z-x)\ dx & 0 & 1/(z-y)
\end{vmatrix},
\end{multline}
where $\kappa = \det H_{\vecn,\vecm}$ is the determinant of the top left block
of \eqref{proof3:startingmatrix0bis}.

Now we expand the integral in the right hand side of
\eqref{proof3:startingmatrix1}.
By definition, each $H_{n_k,m_l}^{(k,l)}$ is a Hankel matrix whose $(i,j)$th
entry is the $(i+j)$th moment with respect to the weight function $W_{k,l}(x) =
w_{1,k}(x)w_{2,l}(x)$, for $k,l=1,2$:
\begin{equation}\label{proof3:defmoments}
(H_{n_k,m_l}^{(k,l)})_{i,j} = \int x^{i+j} W_{k,l}(x)\ dx.
\end{equation}
Note that the integration variable $x$ in \eqref{proof3:defmoments} is just a
\lq dummy variable\rq\ which we can give any name we want. By choosing the same
integration variable $x_1$ for all entries in the first column of
\eqref{proof3:startingmatrix1}, we can rewrite \eqref{proof3:startingmatrix1}
as follows:
\begin{equation*}\label{proof3:matrix1}
=\frac{1}{\kappa}\left|\begin{array}{cccc|ccc|cc} \int W_{1,1}(x_1)\ dx_1 & * &
\ldots &
* &
* &
\ldots & * & 1 & 0\\
\vdots & \vdots & & \vdots & \vdots & & \vdots & \vdots & \vdots\\
\int x_1^{n_1-1} W_{1,1}(x_1)\ dx_1 & * & \ldots & * & * & \ldots & * &
y^{n_1-1} &
0 \\
\hline \int W_{2,1}(x_1)\ dx_1 & * & \ldots & * & * &
\ldots & * & 0 & 1 \\
\vdots & \vdots & & \vdots & \vdots & & \vdots & \vdots & \vdots\\
\int x_1^{n_2-1} W_{2,1}(x_1)\ dx_1 & * & \ldots & * & * & \ldots & * & 0 &
y^{n_2-1}\\
\hline \int \frac{W_{1,1}(x_1)}{z-x_1}\ dx_1 & * & \ldots & * & * & \ldots & *
& \frac{1}{z-y} &
0\\
\int \frac{W_{2,1}(x_1)}{z-x_1}\ dx_1 & * & \ldots & * & * & \ldots & * & 0 &
\frac{1}{z-y}
\end{array}\right|
\end{equation*}
where the entries denoted with $*$ are not important at this stage of the
proof. Using the multi-linearity of the determinant with respect to its first
column, we can take out the integration with respect to $x_1$:
\begin{equation*}\label{proof3:matrix2}
=\frac{1}{\kappa}\int\ \left|\begin{array}{cccc|ccc|cc} w_{1,1}(x_1) &
* & \ldots &
* &
* &
\ldots & * & 1 & 0\\
\vdots & \vdots & & \vdots & \vdots & & \vdots & \vdots & \vdots\\
x_1^{n_1-1}w_{1,1}(x_1) & * & \ldots & * & * & \ldots &
* & y^{n_1-1} &
0 \\
\hline w_{1,2}(x_1) & * & \ldots & * & * &
\ldots & * & 0 & 1 \\
\vdots & \vdots & & \vdots & \vdots & & \vdots & \vdots & \vdots\\
x_1^{n_2-1}w_{1,2}(x_1) & * & \ldots & * & * & \ldots & * & 0 &
y^{n_2-1}\\
\hline \frac{w_{1,1}(x_1)}{z-x_1} & * & \ldots & * & * & \ldots & * &
\frac{1}{z-y} &
0\\
\frac{w_{1,2}(x_1)}{z-x_1} & * & \ldots & * & * & \ldots & * & 0 &
\frac{1}{z-y}
\end{array}\right|\ w_{2,1}(x_1)dx_1.
\end{equation*}
Here we also used that $W_{1,1}(x) = w_{1,1}(x)w_{2,1}(x)$ and $W_{2,1}(x) =
w_{1,2}(x)w_{2,1}(x)$, allowing us to extract the common factor $w_{2,1}(x_1)$
from the first column. Applying the same technique to the other columns of the
first block column we obtain
\begin{equation*}\label{proof3:matrix3}
=\frac{1}{\kappa}\int\ldots\int\ \left|\begin{array}{ccc|ccc|cc}
w_{1,1}(x_1)  & \ldots & w_{1,1}(x_{m_1}) & * & \ldots & * & 1 & 0\\
\vdots & & \vdots & \vdots & & \vdots & \vdots & \vdots\\
x_1^{n_1-1}w_{1,1}(x_1)  & \ldots & x_{m_1}^{n_1-1}w_{1,1}(x_{m_1}) &
* & \ldots &
* & y^{n_1-1} &
0 \\
\hline w_{1,2}(x_1) & \ldots & w_{1,2}(x_{m_1}) & * &
\ldots & * & 0 & 1 \\
\vdots & & \vdots & \vdots & & \vdots & \vdots & \vdots\\
x_1^{n_2-1}w_{1,2}(x_1) & \ldots & x_{m_1}^{n_2-1}w_{1,2}(x_{m_1}) & * & \ldots
& * & 0 & y^{n_2-1}\\
\hline \frac{w_{1,1}(x_1)}{z-x_1} & \ldots & \frac{w_{1,1}(x_{m_1})}{z-x_{m_1}}
&
* & \ldots & * & \frac{1}{z-y} &
0\\
\frac{w_{1,2}(x_1)}{z-x_1} & \ldots & \frac{w_{1,2}(x_{m_1})}{z-x_{m_1}} & * &
\ldots &
* & 0 & \frac{1}{z-y}
\end{array}\right|\ \prod_{j=1}^{m_1}\left(x_j^{j-1}w_{2,1}(x_j)dx_j\right).
\end{equation*}
To get rid of the factor $\prod_{j=1}^{m_1}x_j^{j-1}$ one can apply the
following well-known symmetrization trick. Sum the above expression over all
permutations $x_{\sigma_1},\ldots,x_{\sigma_{m_1}}$ of $x_1,\ldots,x_{m_1}$.
For each term in this sum, apply a column permutation to arrange the columns of
the determinant back in their original form; this releases a factor
$(-1)^{\sigma}$. But then we get the expansion of a Vandermonde determinant:
$$ \sum_{\sigma} (-1)^{\sigma}\prod_{j=1}^{m_1}x_{\sigma_j}^{j-1} = \prod_{1\leq i<j\leq m_1}
(x_j-x_i).
$$
So we get
\begin{multline*}
=\frac{1}{\kappa(m_1)!}\int\ldots\int\ \left|\begin{array}{ccc|ccc|cc}
w_{1,1}(x_1)  & \ldots & w_{1,1}(x_{m_1}) & * & \ldots & * & 1 & 0\\
\vdots & & \vdots & \vdots & & \vdots & \vdots & \vdots\\
x_1^{n_1-1}w_{1,1}(x_1)  & \ldots & x_{m_1}^{n_1-1}w_{1,1}(x_{m_1}) &
* & \ldots &
* & y^{n_1-1} &
0 \\
\hline w_{1,2}(x_1) & \ldots & w_{1,2}(x_{m_1}) & * &
\ldots & * & 0 & 1 \\
\vdots & & \vdots & \vdots & & \vdots & \vdots & \vdots\\
x_1^{n_2-1}w_{1,2}(x_1) & \ldots & x_{m_1}^{n_2-1}w_{1,2}(x_{m_1}) & * & \ldots
& * & 0 & y^{n_2-1}\\
\hline \frac{w_{1,1}(x_1)}{z-x_1} & \ldots & \frac{w_{1,1}(x_{m_1})}{z-x_{m_1}}
&
* & \ldots & * & \frac{1}{z-y} &
0\\
\frac{w_{1,2}(x_1)}{z-x_1} & \ldots & \frac{w_{1,2}(x_{m_1})}{z-x_{m_1}} & * &
\ldots &
* & 0 & \frac{1}{z-y}
\end{array}\right|
\\ \prod_{1\leq i<j\leq m_1} (x_j-x_i)
\prod_{j=1}^{m_1} w_{2,1}(x_j) \prod_{j=1}^{m_1} dx_j.
\end{multline*}
Applying the same techniques to the second block column, we find
\begin{multline}\label{proof3:beforeclaim}
=\frac{1}{\kappa (m_1)!(m_2)!}\int\ldots\int\ \left|\begin{array}{ccc|cc}
w_{1,1}(x_1) &
\ldots & w_{1,1}(x_{n}) & 1 & 0 \\
\vdots & & \vdots & \vdots & \vdots \\
x_1^{n_1-1} w_{1,1}(x_1) & \ldots & x_{n}^{n_1-1} w_{1,1}(x_{n}) & y^{n_1-1} & 0 \\
\hline w_{1,2}(x_1) & \ldots & w_{1,2}(x_{n}) & 0 & 1 \\
\vdots & & \vdots & \vdots & \vdots \\
x_1^{n_2-1} w_{1,2}(x_1) & \ldots & x_{n}^{n_2-1} w_{1,2}(x_{n}) & 0 & y^{n_2-1} \\
\hline \frac{w_{1,1}(x_1)}{z-x_1} & \ldots & \frac{w_{1,1}(x_{n})}{z-x_{n}} & \frac{1}{z-y} & 0 \\
\frac{w_{1,2}(x_1)}{z-x_1} & \ldots & \frac{w_{1,2}(x_{n})}{z-x_{n}} & 0 &
\frac{1}{z-y}
\end{array}\right|\\ \prod_{1\leq i<j\leq m_1} (x_j-x_i)
\prod_{m_1+1\leq i<j\leq n} (x_j-x_i) \prod_{j=1}^{m_1}w_{2,1}(x_j)
\prod_{j=m_1+1}^{n} w_{2,2}(x_j)\prod_{j=1}^n dx_j.
\end{multline}

We claim that this can be rewritten as
\begin{multline}\label{proof3:righthandside:reduced}
=\frac{1}{(z-y)^2}\frac{1}{\kappa(m_1)!(m_2)!}\int\ldots\int\
\frac{\prod_{j=1}^n (y-x_j)}{\prod_{j=1}^n (z-x_j)}\ \left|\begin{array}{ccc}
w_{1,1}(x_1) &
\ldots & w_{1,1}(x_{n}) \\
\vdots & & \vdots  \\
x_1^{n_1-1} w_{1,1}(x_1) & \ldots & x_{n}^{n_1-1} w_{1,1}(x_{n}) \\
\hline w_{1,2}(x_1) & \ldots & w_{1,2}(x_{n}) \\
\vdots & & \vdots \\
x_1^{n_2-1} w_{1,2}(x_1) & \ldots & x_{n}^{n_2-1} w_{1,2}(x_{n})
\end{array}\right|\\ \prod_{1\leq i<j\leq m_1} (x_j-x_i) \prod_{m_1+1\leq
i<j\leq n} (x_j-x_i) \prod_{j=1}^{m_1}w_{2,1}(x_j) \prod_{j=m_1+1}^{n}
w_{2,2}(x_j)\prod_{j=1}^n dx_j.
\end{multline}
This follows from standard determinant manipulations and the exact formula for
a Cauchy-Vandermonde determinant; see the appendix.

Now up to a constant factor times $1/(z-y)^2$,
\eqref{proof3:righthandside:reduced} equals
\eqref{proof3:lefthandside:reduced}. This proves \eqref{mainR3a}, up to a
constant factor; but this constant factor must be 1 as can be seen by matching
the leading terms on both sides of \eqref{mainR3a}.

\begin{remark} By matching the leading constants in
\eqref{proof3:lefthandside:reduced} and \eqref{proof3:righthandside:reduced},
we obtain the following expression for the normalization constant $Z_n$ in
\eqref{determinantalpointprocess}:
$$ Z_n = \kappa n!
$$
where we recall that $\kappa$ is the block Hankel determinant $\kappa = \det
H_{\vecn,\vecm}$. This identity also follows directly from a \lq generalized
Cauchy-Binet identity\rq, see e.g. \cite{Kui}.
\end{remark}

\subsection{Proof of Theorems \ref{theorem:main1} and \ref{theorem:main2}}
\label{subsection:prooftheorem:1and2}

Having established Theorem \ref{theorem:main3}, Theorems \ref{theorem:main1}
and \ref{theorem:main2} now follow as simple limiting cases. For Theorem
\ref{theorem:main1} we let $z\to\infty$ in \eqref{mainR3a} and observe from
\eqref{def:Rn} that
\begin{multline*}
\det R_n(z,y) = \det\left(
\begin{pmatrix}I_p & 0 \end{pmatrix}Y^{-1}(z)Y(y)\begin{pmatrix}I_p \\
0
\end{pmatrix}\right) = \det\left(
\begin{pmatrix}(Y^{-1})_{1,1}(z) & (Y^{-1})_{1,2}(z) \end{pmatrix}\begin{pmatrix}Y_{1,1}(y) \\
Y_{2,1}(y)
\end{pmatrix}\right) \\ = z^{-n} \det Y_{1,1}(y)+O(z^{-n-1}),
\end{multline*}
where in the last step we used the Cauchy-Binet identity and the asymptotics
\eqref{asymptoticconditionY0}. This establishes Theorem \ref{theorem:main1}.
Theorem \ref{theorem:main2} can be obtained in a similar way by letting
$y\to\infty$ in \eqref{mainR3b}.

\subsection{Proof of Theorem \ref{theorem:main4}(c)}
\label{subsection:prooftheorem:ratio:more}

The proof that we have given for Theorem \ref{theorem:main3} can be extended to
prove Theorem \ref{theorem:main4}(c) as well. The key observation is that the
matrix in the right hand side of \eqref{main4c} can be written in Schur
complement form as well. For example, for $K=L=2$ and $p=q=2$ one has that
\begin{equation*}
\begin{pmatrix}
\frac{1}{z_1-y_1} R_n(z_1,y_1) & \frac{1}{z_1-y_2} R_n(z_1,y_2) \\
\frac{1}{z_2-y_1} R_n(z_2,y_1) & \frac{1}{z_2-y_2} R_n(z_2,y_2)
\end{pmatrix}
\end{equation*}
is the Schur complement of the matrix.
\begin{equation}\label{proof4:startingmatrix}
\left(\begin{array}{cc|cccc} H_{n_1,m_1}^{(1,1)} & H_{n_1,m_2}^{(1,2)} &
(\mathbf{y}_1)_{n_1} & \mathbf{0} &
(\mathbf{y}_2)_{n_1} & \mathbf{0}  \\
H_{n_2,m_1}^{(2,1)} & H_{n_2,m_2}^{(2,2)} & \mathbf{0} & (\mathbf{y_1})_{n_2} &
\mathbf{0} & (\mathbf{y_2})_{n_2} \\
\hline \vphantom{\overset{O}{\sim}}\int\mathbf{x}^T_{m_1}W_{1,1}(x)/(z_1-x)\ dx & \int\mathbf{x}^T_{m_2}W_{1,2}(x)/(z_1-x)\ dx & 1/(z_1-y_1) & 0 & 1/(z_1-y_2) & 0\\
\int\mathbf{x}^T_{m_1}W_{2,1}(x)/(z_1-x)\ dx &
\int\mathbf{x}^T_{m_2}W_{2,2}(x)/(z_1-x)\ dx & 0 & 1/(z_1-y_1) & 0 &
1/(z_1-y_2)\\
\int\mathbf{x}^T_{m_1}W_{1,1}(x)/(z_2-x)\ dx & \int\mathbf{x}^T_{m_2}W_{1,2}(x)/(z_2-x)\ dx & 1/(z_2-y_1) & 0 & 1/(z_2-y_2) & 0\\
\int\mathbf{x}^T_{m_1}W_{2,1}(x)/(z_2-x)\ dx &
\int\mathbf{x}^T_{m_2}W_{2,2}(x)/(z_2-x)\ dx & 0 & 1/(z_2-y_1) & 0 &
1/(z_2-y_2)
\end{array}\right).
\end{equation}
Compare with \eqref{proof3:startingmatrix0bis}.

The determinant of this Schur complement can be manipulated in the same way as
in Section~\ref{subsection:prooftheorem:ratio}. This leads to the following
analogue of \eqref{proof3:beforeclaim}:
\begin{multline}\label{proof4:beforeclaim}
\frac{1}{\kappa (m_1)!(m_2)!}\int\ldots\int\ \left|\begin{array}{ccc|cccc}
w_{1,1}(x_1) &
\ldots & w_{1,1}(x_{n}) & 1 & 0 & 1 & 0 \\
\vdots & & \vdots & \vdots & \vdots & \vdots & \vdots \\
x_1^{n_1-1} w_{1,1}(x_1) & \ldots & x_{n}^{n_1-1} w_{1,1}(x_{n}) & y_1^{n_1-1} & 0 & y_2^{n_1-1} & 0 \\
\hline w_{1,2}(x_1) & \ldots & w_{1,2}(x_{n}) & 0 & 1 & 0 & 1 \\
\vdots & & \vdots & \vdots & \vdots & \vdots & \vdots \\
x_1^{n_2-1} w_{1,2}(x_1) & \ldots & x_{n}^{n_2-1} w_{1,2}(x_{n}) & 0 & y_1^{n_2-1} & 0 & y_2^{n_2-1} \\
\hline \frac{w_{1,1}(x_1)}{z_1-x_1} & \ldots & \frac{w_{1,1}(x_{n})}{z_1-x_{n}} & \frac{1}{z_1-y_1} & 0 & \frac{1}{z_1-y_2} & 0 \\
\frac{w_{1,2}(x_1)}{z_1-x_1} & \ldots & \frac{w_{1,2}(x_{n})}{z_1-x_{n}} & 0 &
\frac{1}{z_1-y_1} & 0 & \frac{1}{z_1-y_2} \\
\frac{w_{1,1}(x_1)}{z_2-x_1} & \ldots & \frac{w_{1,1}(x_{n})}{z_2-x_{n}} & \frac{1}{z_2-y_1} & 0 & \frac{1}{z_2-y_2} & 0 \\
\frac{w_{1,2}(x_1)}{z_2-x_1} & \ldots & \frac{w_{1,2}(x_{n})}{z_2-x_{n}} & 0 &
\frac{1}{z_2-y_1} & 0 & \frac{1}{z_2-y_2}
\end{array}\right|\\ \prod_{1\leq i<j\leq m_1} (x_j-x_i)
\prod_{m_1+1\leq i<j\leq n} (x_j-x_i) \prod_{j=1}^{m_1}w_{2,1}(x_j)
\prod_{j=m_1+1}^{n} w_{2,2}(x_j)\prod_{j=1}^n dx_j,
\end{multline}
and the following analogue of \eqref{proof3:righthandside:reduced}:
\begin{multline}\label{proof4:righthandside:reduced}
\left(\frac{(z_1-z_2)(y_2-y_1)}{\prod_{i,j=1}^{2}
(z_i-y_j)}\right)^2\frac{1}{\kappa (m_1)!(m_2)!}\int\ldots\int\
\frac{\prod_{i=1}^2\prod_{j=1}^n (y_i-x_j)}{\prod_{i=1}^2\prod_{j=1}^n
(z_i-x_j)}\left|\begin{array}{ccc} w_{1,1}(x_1) &
\ldots & w_{1,1}(x_{n}) \\
\vdots & & \vdots  \\
x_1^{n_1-1} w_{1,1}(x_1) & \ldots & x_{n}^{n_1-1} w_{1,1}(x_{n}) \\
\hline w_{1,2}(x_1) & \ldots & w_{1,2}(x_{n}) \\
\vdots & & \vdots \\
x_1^{n_2-1} w_{1,2}(x_1) & \ldots & x_{n}^{n_2-1} w_{1,2}(x_{n})
\end{array}\right|\\ \prod_{1\leq i<j\leq m_1} (x_j-x_i)
\prod_{m_1+1\leq i<j\leq n} (x_j-x_i) \prod_{j=1}^{m_1}w_{2,1}(x_j)
\prod_{j=m_1+1}^{n} w_{2,2}(x_j)\prod_{j=1}^n dx_j.
\end{multline}
Here the transition from \eqref{proof4:beforeclaim} to
\eqref{proof4:righthandside:reduced} follows again from standard determinant
manipulations and the exact formula for a Cauchy-Vandermonde determinant, see
the appendix.

Formula \eqref{proof4:righthandside:reduced} gives an expression for the
determinant in the right hand side of \eqref{main4c}. Here the leftmost factor
\begin{equation}\label{proof4:prefactor} \left(\frac{(z_1-z_2)(y_2-y_1)}{\prod_{i,j=1}^{2} (z_i-y_j)}\right)^2\end{equation} is
canceled by multiplying with the first factor in the right hand side in
\eqref{main4c}.  Finally, by also expanding the left hand side of
\eqref{main4c} in the same way as before, cf.\
\eqref{proof3:lefthandside:reduced}, Theorem \ref{theorem:main4}(c) follows.

\section{Proofs, part 2}
\label{section:proofmaintheorem:ratios}

In this section we prove Theorems \ref{theorem:main4}(a)--(b) and
\ref{theorem:main5} about averages of general products and ratios of
characteristic polynomials. To this end we use a similar strategy as in
Baik-Deift-Strahov~\cite{BDS}, i.e., we make use of the already established
Theorems \ref{theorem:main1}--\ref{theorem:main3}, together with adaptations of
the classical results by Christoffel and Uvarov. We mention that these
Christoffel-Uvarov type results will be valid for arbitrary vector orthogonal
polynomials, i.e., the weight matrix $W(x)$ does not need to have the rank-one
factorization \eqref{weightmatrixrankone}.

Throughout this section we assume that $|\vecn| = |\vecm|$ and we let $W(x)$ be
an arbitrary weight matrix of size $p\times q$. We let $Y(z)$ be the solution
of the corresponding RH problem and $K_n(x,y)$ the kernel
\eqref{def:kernel:block}. We also consider the modified weight matrix
\begin{equation}\label{mixedchristoffel:Wtilde}\tilde W(x) = \frac{\prod_{k=1}^{K} (x-y_k)}{\prod_{l=1}^{L} (x-z_l)}W(x),\end{equation}
and we denote the corresponding RH matrix and kernel with $\tilde{Y}(z)$ and
$\tilde{K}_n(x,y)$. Here we assume that $y_1,\ldots,y_K\in\cee$,
$z_1,\ldots,z_L\in\cee\setminus\er$ and the numbers in the set
$(y_1,\ldots,y_K$, $z_1,\ldots,z_L)$ are pairwise distinct. Our aim will be to
relate the matrices $\tilde{Y}$ and $Y$, and similarly for the kernels
$\tilde{K}_n$ and $K_n$.

Although Theorem \ref{theorem:main4} is just a special case of Theorem
\ref{theorem:main5}, for reasons of comprehension we discuss the former first,
see Sections \ref{subsection:prooftheorem:4a} and
\ref{subsection:prooftheorem:4b}. The proof of Theorem \ref{theorem:main5} is
then only briefly sketched in Section \ref{subsection:prooftheorem:5}.

\subsection{Proof of Theorem \ref{theorem:main4}(a)}
\label{subsection:prooftheorem:4a}

Recall that \emph{Christoffel's formula} expresses the orthogonal polynomials
with respect to a polynomially modified weight function $\tilde w(x) =
\prod_{k=1}^{K} (x-y_k)w(x)$ in terms of those with respect to the original
weight function $w(x)$ \cite{BDS,Sz}. We now state a generalization of this
formula to the context of vector orthogonal polynomials. Throughout, we will
use the notations in Definition~\ref{def:chainmultiindices}.

\begin{proposition}\label{proposition:christoffel}
Assume that $L=0$ in \eqref{mixedchristoffel:Wtilde}, so that
\begin{equation}\label{christoffel:Wtilde}\tilde W(x) = \prod_{k=1}^{K} (x-y_k)W(x).\end{equation}
 Then $\tilde{Y}_{1,1}(y)$ is the Schur
complement of the matrix
\begin{equation}\label{christoffel}
\frac{1}{\prod_{k=1}^{K}(y-y_k)} \left(\begin{array}{ccc|c}
(Y_{1,1})_{\vecn_0,\vecm_0}(y_1) & \ldots & (Y_{1,1})_{\vecn_0,\vecm_0}(y_{K}) & (Y_{1,1})_{\vecn_0,\vecm_0}(y) \\
\vdots & & \vdots & \vdots \\
(Y_{1,1})_{\vecn_{K-1},\vecm_{K-1}}(y_1) & \ldots &
(Y_{1,1})_{\vecn_{K-1},\vecm_{K-1}}(y_K) & (Y_{1,1})_{\vecn_{K-1},\vecm_{K-1}}(y) \\
\hline (Y_{1,1})_{\vecn_{K},\vecm_{K}}(y_1) & \ldots &
(Y_{1,1})_{\vecn_{K},\vecm_{K}}(y_{K}) & (Y_{1,1})_{\vecn_{K},\vecm_{K}}(y)
\end{array}\right)
\end{equation}
with respect to its bottom right $p$ by $p$ submatrix.
\end{proposition}

\bewijs. Denote with $S(y)$ the Schur complement of the matrix
\eqref{christoffel}. It is clear that $S(y)$ is a matrix of size $p$ by $p$. We
claim that the entries of $S(y)$ are polynomials in $y$, and moreover that
\begin{equation}\label{proof:christoffel1}
S(y) = (I_{p}+O(1/y))\diag(y^{n_1},\ldots,y^{n_p}),\qquad y\to\infty.
\end{equation}
To prove these claims, first note that all entries in the last block column of
\eqref{christoffel} are polynomials in $y$. Now we use the expression
\eqref{Schur:entrywise3} for the entries of the Schur complement of the matrix
\eqref{defSchurcomp1}. It is clear that if we apply this formula to the input
matrix \eqref{christoffel}, then the determinant in the numerator of
\eqref{Schur:entrywise3} is zero whenever $y=y_k$ for certain $k=1,\ldots,K$.
This means that all entries of $S(y)$ are polynomials divisible by
$\prod_{k=1}^{K} (y-y_k)$, which on multiplying with the prefactor
$1/\prod_{k=1}^{K} (y-y_k)$ in \eqref{christoffel} indeed shows that $S(y)$ is
a polynomial matrix.

To establish \eqref{proof:christoffel1}, observe that the leading term of
$S(y)$ comes from the bottom right $p$ by $p$ block in \eqref{christoffel},
which has the asymptotics
\begin{equation*}\label{proof:christoffel2}
(Y_{1,1})_{\vecn_{K},\vecm_{K}}(y) =
(I_{p}+O(1/y))\diag(y^{n_1+K},\ldots,y^{n_p+K}),\qquad y\to\infty.
\end{equation*}
Thus on dividing by the prefactor $1/\prod_{k=1}^{K} (y-y_k)$ in
\eqref{christoffel}, we obtain \eqref{proof:christoffel1}.

Next, we check the orthogonality conditions
\begin{equation}\label{proof:christoffel3}
\int S(x)\tilde{W}(x)\vecQ(x)\ dx = \mathbf{0}
\end{equation}
for all $\vecQ\in\Pee_{\vecm}$, where we recall the definition of $\tilde W(x)$
in \eqref{christoffel:Wtilde}. To prove this, note that the factor
$\prod_{k=1}^{K} (x-y_k)$ in \eqref{christoffel:Wtilde} cancels with the
prefactor in \eqref{christoffel}, and by the multi-linearity of determinants
the remaining integral can be taken into the last block column of
\eqref{christoffel}. So the left hand side of \eqref{proof:christoffel3} is the
Schur complement of the matrix
\begin{equation}\label{proof:christoffel4}
\left(\begin{array}{ccc|c}
(Y_{1,1})_{\vecn_0,\vecm_0}(y_1) & \ldots & (Y_{1,1})_{\vecn_0,\vecm_0}(y_{K}) & \int (Y_{1,1})_{\vecn_0,\vecm_0}(x)W(x)\vecQ(x)\ dx \\
\vdots & & \vdots & \vdots \\
(Y_{1,1})_{\vecn_{K-1},\vecm_{K-1}}(y_1) & \ldots &
(Y_{1,1})_{\vecn_{K-1},\vecm_{K-1}}(y_K) & \int (Y_{1,1})_{\vecn_{K-1},\vecm_{K-1}}(x)W(x)\vecQ(x)\ dx \\
\hline (Y_{1,1})_{\vecn_{K},\vecm_{K}}(y_1) & \ldots &
(Y_{1,1})_{\vecn_{K},\vecm_{K}}(y_{K}) & \int
(Y_{1,1})_{\vecn_{K},\vecm_{K}}(x)W(x)\vecQ(x)\ dx
\end{array}\right).
\end{equation}
But each of the matrices $Y_{1,1}$ in the last column of
\eqref{proof:christoffel4} satisfies the orthogonality relations with respect
to $\vecQ\in\Pee_{\vecm}$. So the last column of \eqref{proof:christoffel4} is
zero and therefore the Schur complement of this matrix is zero as well. This
establishes \eqref{proof:christoffel3}.

The proof of the proposition now follows from \eqref{proof:christoffel1} and
\eqref{proof:christoffel3}, see \eqref{MOPtype2}. $\bol$

Having established Theorem \ref{theorem:main1} and Proposition
\ref{proposition:christoffel}, and taking into account the property
\eqref{Schurdeterminants}, Theorem \ref{theorem:main4}(a) can now be proved by
induction on $K$, in exactly the same way as in \cite[Proof of Th.~4.3]{BDS}.
We refer to the latter reference for details.

\subsection{Proof of Theorem \ref{theorem:main4}(b)}
\label{subsection:prooftheorem:4b}

Recall that \emph{Uvarov's formula} expresses the orthogonal polynomials with
respect to an inverse-polynomially modified weight function $\tilde w(x) =
w(x)/\prod_{l=1}^{L} (x-z_l)$ in terms of those with respect to the original
weight function $w(x)$ \cite{BDS,Uv}. We now state a generalization of this
formula to the context of vector orthogonal polynomials.

\begin{proposition}\label{proposition:uvarov}
Assume that $K=0$ in \eqref{mixedchristoffel:Wtilde}, so that
\begin{equation}\label{uvarov:Wtilde}\tilde W(x) = \frac{1}{\prod_{l=1}^{L} (x-z_l)}W(x).\end{equation}
Then $\tilde{Y}_{2,1}(z)$ is the Schur complement of the matrix
\begin{equation}\label{uvarov:Y21}
\left(\begin{array}{ccc|c}
(Y_{2,2})_{\vecn_0,\vecm_0}(z_1) & \ldots & (Y_{2,2})_{\vecn_0,\vecm_0}(z_{L}) & (Y_{2,1})_{\vecn_0,\vecm_0}(z) \\
\vdots & & \vdots & \vdots \\
(Y_{2,2})_{\vecn_{-L+1},\vecm_{-L+1}}(z_1) & \ldots &
(Y_{2,2})_{\vecn_{-L+1},\vecm_{-L+1}}(z_L) & (Y_{2,1})_{\vecn_{-L+1},\vecm_{-L+1}}(z) \\
\hline (Y_{2,2})_{\vecn_{-L},\vecm_{-L}}(z_{1}) & \ldots &
(Y_{2,2})_{\vecn_{-L},\vecm_{-L}}(z_{L}) & (Y_{2,1})_{\vecn_{-L},\vecm_{-L}}(z)
\end{array}\right)
\end{equation}
with respect to its bottom right $q$ by $p$ submatrix, and $\tilde{Y}_{2,2}(z)$
is the Schur complement of the matrix
\begin{equation}\label{uvarov:Y22}
\frac{1}{\prod_{l=1}^{L} (z-z_l)}\left(\begin{array}{ccc|c}
(Y_{2,2})_{\vecn_0,\vecm_0}(z_1) & \ldots & (Y_{2,2})_{\vecn_0,\vecm_0}(z_{L}) & (Y_{2,2})_{\vecn_0,\vecm_0}(z) \\
\vdots & & \vdots & \vdots \\
(Y_{2,2})_{\vecn_{-L+1},\vecm_{-L+1}}(z_1) & \ldots &
(Y_{2,2})_{\vecn_{-L+1},\vecm_{-L+1}}(z_L) & (Y_{2,2})_{\vecn_{-L+1},\vecm_{-L+1}}(z) \\
\hline (Y_{2,2})_{\vecn_{-L},\vecm_{-L}}(z_{1}) & \ldots &
(Y_{2,2})_{\vecn_{-L},\vecm_{-L}}(z_{L}) & (Y_{2,2})_{\vecn_{-L},\vecm_{-L}}(z)
\end{array}\right)
\end{equation}
with respect to its bottom right $q$ by $q$ submatrix.
\end{proposition}

\bewijs. Denote with $S(z)$ the Schur complement of the matrix
\eqref{uvarov:Y21}. It is clear that $S(z)$ is a polynomial matrix of size $q$
by $p$, which has the required degree structure
\begin{equation}\label{proof:uvarov1}
S(z) = (O(1/z))\diag(z^{n_1},\ldots,z^{n_p}),\qquad z\to\infty.
\end{equation}
On the other hand, denote with $T(z)$ the Schur complement of the matrix
\eqref{uvarov:Y22}. From the asymptotic condition
\begin{equation*}\label{proof:uvarov2}
(Y_{2,2})_{\vecn_{-L},\vecm_{-L}}(z) =
(I_{q}+O(1/z))\diag(z^{L-m_1},\ldots,z^{L-m_q}),\qquad z\to\infty,
\end{equation*}
it is clear that $T(z)$ has the required asymptotics
\begin{equation}\label{proof:uvarov3}
T(z) = (I_{q}+O(1/z))\diag(z^{-m_1},\ldots,z^{-m_q}),\qquad z\to\infty.
\end{equation}
Taking into account \eqref{proof:uvarov1} and \eqref{proof:uvarov3}, the
proposition will follow if we can show that
\begin{equation}\label{proof:uvarov4}
T(z) = -\frac{1}{2\pi i}\int \frac{S(x)\tilde W(x)}{z-x}\ dx.
\end{equation}
Consider the left hand side of \eqref{proof:uvarov4}. From the formula
$Y_{2,2}(z) = -\frac{1}{2\pi i}\int \frac{Y_{2,1}(x)W(x)}{z-x}\ dx$ we see that
all entries of matrix \eqref{uvarov:Y22} are integrals. Choose a common
integration variable $x_i$ for all entries in the $i$th row. Consider the
coefficients in the partial fraction decomposition
\begin{equation}\label{uvarov:pfd} \sum_{l=1}^L \frac{c_l}{x-z_l}+\frac{1}{x-z}
= \frac{c}{(x-z)\prod_{l=1}^L (x-z_l)},
\end{equation} where $c_l, c$ are functions of $\{z_1,\ldots,z_L,z\}$ but not of $x$.
Consider the column operation where to the last
block column of matrix \eqref{uvarov:Y22} we add $c_l$ times the $l$th block
column, $l=1,\ldots,L$. Performing these operations inside the integrals, this
causes the integrand in the last block column to be multiplied with
$c/\prod_{l=1}^{L} (x-z_l) = \prod_{l=1}^{L} (z-z_l)/\prod_{l=1}^{L} (x-z_l)$.
Here the factor $\prod_{l=1}^{L} (z-z_l)$ can be taken out of the integrand,
canceling the prefactor in \eqref{uvarov:Y22}, while the factor
$W(x)/\prod_{l=1}^{L} (x-z_l)$ in the integrand is nothing but $\tilde{W}(x)$.
So we obtain precisely the right hand side of \eqref{proof:uvarov4}. This
establishes \eqref{proof:uvarov4}, and thereby the proposition is proved.
 $\bol$

Having established Theorem \ref{theorem:main2} and Proposition
\ref{proposition:uvarov}, and taking into account the property
\eqref{Schurdeterminants}, Theorem \ref{theorem:main4}(b) can now be proved by
induction on $L$, in exactly the same way as in \cite[Proof of Th.~4.10]{BDS}.
We refer to the latter reference for details.

\subsection{Proof of Theorem \ref{theorem:main5}}
\label{subsection:prooftheorem:5}

Finally we prove Theorem \ref{theorem:main5}. The key to the proof of Theorem
\ref{theorem:main5}(a) is the following result.

\begin{proposition}\label{proposition:mixedchristoffel}
Assume that $K\geq L$ in \eqref{mixedchristoffel:Wtilde}. Then
$\tilde{Y}_{1,1}(y)$ is the Schur complement of the matrix
\begin{equation}\label{mixedchristoffel}
\frac{\prod_{l=1}^{L}(y-z_l)}{\prod_{k=1}^{K}(y-y_k)}
\left(\begin{array}{ccc|c}
\frac{1}{z_1-y_1} R_n(z_1,y_1) & \ldots & \frac{1}{z_1-y_K} R_n(z_1,y_K) & \frac{1}{z_1-y} R_n(z_1,y) \\
\vdots & & \vdots & \vdots \\
\frac{1}{z_L-y_1} R_n(z_L,y_1) & \ldots & \frac{1}{z_L-y_K} R_n(z_L,y_K) & \frac{1}{z_L-y} R_n(z_L,y) \\
(Y_{1,1})_{\vecn_{0},\vecm_{0}}(y_1) & \ldots &
(Y_{1,1})_{\vecn_{0},\vecm_{0}}(y_{K}) & (Y_{1,1})_{\vecn_{0},\vecm_{0}}(y)\\
\vdots & & \vdots & \vdots \\ \hline (Y_{1,1})_{\vecn_{K-L},\vecm_{K-L}}(y_1) &
\ldots & (Y_{1,1})_{\vecn_{K-L},\vecm_{K-L}}(y_{K}) &
(Y_{1,1})_{\vecn_{K-L},\vecm_{K-L}}(y)
\end{array}\right)
\end{equation}
with respect to its bottom right $p$ by $p$ submatrix. Here all $R_n$ matrices
are taken with respect to the pair of multi-indices $\vecn_{0},\vecm_{0}$.
\end{proposition}

The proof follows the same lines as the proof of Proposition
\ref{proposition:christoffel}, and we omit the details. We suffice to mention
that the proof makes use of the vanishing property for the matrix $R_n$ in
Proposition \ref{proposition:vanishingprop:dual}.

From Proposition \ref{proposition:mixedchristoffel}, Theorem
\ref{theorem:main5}(a) can now be proved by induction on $K$, $K\geq L$ (for
$L$ fixed), with induction basis $K=L$ corresponding to Theorem
\ref{theorem:main4}(c).

Finally, we note that the proof of Theorem \ref{theorem:main5}(b) hinges on the
following fact.

\begin{proposition}\label{proposition:mixeduvarov}
Assume that $L\geq K$ in \eqref{mixedchristoffel:Wtilde}. Then
$\tilde{Y}_{2,1}(z)$ is the Schur complement of the matrix
\begin{equation}\label{mixeduvarov:Y21}
\left(\begin{array}{ccc|c}
\frac{1}{z_1-y_1} L_n(y_1,z_1) & \ldots & \frac{1}{z_L-y_1} L_n(y_1,z_L) & 2\pi i K_n(y_1,z) \\
\vdots & & \vdots & \vdots \\
\frac{1}{z_1-y_K} L_n(y_K,z_1) & \ldots & \frac{1}{z_L-y_K} L_n(y_K,z_L) & 2\pi i K_n(y_K,z) \\
(Y_{2,2})_{\vecn_{0},\vecm_{0}}(z_1) & \ldots &
(Y_{2,2})_{\vecn_{0},\vecm_{0}}(z_{L}) & (Y_{2,1})_{\vecn_{0},\vecm_{0}}(z)\\
\vdots & & \vdots & \vdots \\ \hline (Y_{2,2})_{\vecn_{K-L},\vecm_{K-L}}(z_1) &
\ldots & (Y_{2,2})_{\vecn_{K-L},\vecm_{K-L}}(z_{L}) &
(Y_{2,1})_{\vecn_{K-L},\vecm_{K-L}}(z)
\end{array}\right)
\end{equation}
with respect to its bottom right $q$ by $p$ submatrix, and $\tilde{Y}_{2,2}(z)$
is the Schur complement of the matrix
\begin{equation}\label{mixeduvarov:Y22}
\frac{\prod_{k=1}^{K}(z-y_k)}{\prod_{l=1}^{L}(z-z_l)}
\left(\begin{array}{ccc|c}
\frac{1}{z_1-y_1} L_n(y_1,z_1) & \ldots & \frac{1}{z_L-y_1} L_n(y_1,z_L) & \frac{1}{z-y_1} L_n(y_1,z) \\
\vdots & & \vdots & \vdots \\
\frac{1}{z_1-y_K} L_n(y_K,z_1) & \ldots & \frac{1}{z_L-y_K} L_n(y_K,z_L) & \frac{1}{z-y_K} L_n(y_K,z) \\
(Y_{2,2})_{\vecn_{0},\vecm_{0}}(z_1) & \ldots &
(Y_{2,2})_{\vecn_{0},\vecm_{0}}(z_{L}) & (Y_{2,2})_{\vecn_{0},\vecm_{0}}(z)\\
\vdots & & \vdots & \vdots \\ \hline (Y_{2,2})_{\vecn_{K-L},\vecm_{K-L}}(z_1) &
\ldots & (Y_{2,2})_{\vecn_{K-L},\vecm_{K-L}}(z_{L}) &
(Y_{2,2})_{\vecn_{K-L},\vecm_{K-L}}(z)
\end{array}\right)
\end{equation}
with respect to its bottom right $q$ by $q$ submatrix. Here all $L_n$ matrices
are taken with respect to the pair of multi-indices $\vecn_{0},\vecm_{0}$.
\end{proposition}

\bewijs. The proof follows the same lines as the proof of Proposition
\ref{proposition:uvarov}. The properties \eqref{proof:uvarov1} and
\eqref{proof:uvarov3} are easily extended to the present situation. The main
difficulty lies in proving \eqref{proof:uvarov4}. To this end one can use the
following analogue of the partial fraction decomposition \eqref{uvarov:pfd},
\begin{equation}\label{mixeduvarov:pfd} \sum_{l=1}^L \frac{c_l}{x-z_l}+\frac{1}{x-z}
= \frac{c\prod_{k=1}^K (x-y_k)}{(x-z)\prod_{l=1}^L (x-z_l)},
\end{equation} where $c_l, c$ are functions of $\{y_1,\ldots,y_K,z_1,\ldots,z_L,z\}$
but not of $x$. Observe that these very same coefficients $c_l, c$ also appear
in the more complicated partial fraction decomposition
\begin{equation}\label{mixeduvarov:pfdbis} \sum_{l=1}^L \frac{c_l(x-y_k)}{(z_l-y_k)(x-z_l)}+\frac{x-y_k}{(z-y_k)(x-z)}
= \frac{c\prod_{k=1}^K (x-y_k)}{(x-z)\prod_{l=1}^L (x-z_l)},
\end{equation}
for any $k$. The proof of \eqref{proof:uvarov4} then follows as in the proof of
Proposition \ref{proposition:uvarov}, using the integral representations
$Y_{2,2}(z) = -\frac{1}{2\pi i}\int \frac{Y_{2,1}(x)W(x)}{z-x}\ dx$ (as before)
and
$$ L_n(y,z) = \int \frac{x-y}{x-z}K_n(y,x)W(x)\ dx,$$
and applying the relations \eqref{mixeduvarov:pfd} for the integrals for
$Y_{2,2}$ and \eqref{mixeduvarov:pfdbis} for those for $L_n$. $\bol$

\section*{Appendix: Calculations with Cauchy-Vandermonde determinants}

Define the functions
\begin{equation*}\label{CVM}
f_i(x) = \left\{\begin{array}{ll} x^{i-1},& i=1,\ldots,n, \\
1/(z_{i-n}-x),& i=n+1,\ldots,n+m,
\end{array}\right.
\end{equation*}
where $z_1,\ldots,z_m$ are a given set of numbers. It is well-known that
\begin{equation}\label{cauchy:vdm:det}
\det \begin{pmatrix}f_i(x_j)\end{pmatrix}_{i,j=1}^{n+m} = \frac{\prod_{1\leq i<
j\leq m} (z_i-z_j)\prod_{1\leq i< j\leq n+m} (x_j-x_i)}{\prod_{i,j}(z_i-x_j)}.
\end{equation}
The expression in the left hand side of \eqref{cauchy:vdm:det} is called a
\emph{Cauchy-Vandermonde determinant}.

\subsection*{Proof that \eqref{proof3:beforeclaim} implies
\eqref{proof3:righthandside:reduced}} By applying a suitable row permutation,
the determinant inside \eqref{proof3:beforeclaim} takes the form
\begin{equation*}\label{appendix:1}
(-1)^{n_2}\left|\begin{array}{ccc|cc} w_{1,1}(x_1) &
\ldots & w_{1,1}(x_{n}) & 1 & 0 \\
\vdots & & \vdots & \vdots & \vdots \\
x_1^{n_1-1} w_{1,1}(x_1) & \ldots & x_{n}^{n_1-1} w_{1,1}(x_{n}) & y^{n_1-1} & 0 \\
\frac{w_{1,1}(x_1)}{z-x_1} & \ldots & \frac{w_{1,1}(x_{n})}{z-x_{n}} & \frac{1}{z-y} & 0 \\
\hline w_{1,2}(x_1) & \ldots & w_{1,2}(x_{n}) & 0 & 1 \\
\vdots & & \vdots & \vdots & \vdots \\
x_1^{n_2-1} w_{1,2}(x_1) & \ldots & x_{n}^{n_2-1} w_{1,2}(x_{n}) & 0 & y^{n_2-1} \\
\frac{w_{1,2}(x_1)}{z-x_1} & \ldots & \frac{w_{1,2}(x_{n})}{z-x_{n}} & 0 &
\frac{1}{z-y}
\end{array}\right|.
\end{equation*}
By applying Lagrange expansion with respect to the first block row, this equals
\begin{multline}\label{appendix:Lagrange}
= \sum_{S} (-1)^{S}\begin{vmatrix} w_{1,1}(x_{s_1}) & \ldots & w_{1,1}(x_{s_{n_1}}) & 1 \\
\vdots & & \vdots & \vdots \\
x_{s_{1}}^{n_1-1}w_{1,1}(x_{s_1}) & \ldots &
x_{s_{n_1}}^{n_1-1}w_{1,1}(x_{s_{n_1}}) & y^{n_1-1} \\
\frac{w_{1,1}(x_{s_{1}})}{z-x_{s_1}} & \ldots &
\frac{w_{1,1}(x_{s_{n_1}})}{z-x_{s_{n_1}}}& \frac{1}{z-y}
\end{vmatrix}\\
\begin{vmatrix} w_{1,2}(x_{\bar{s}_1}) & \ldots & w_{1,2}(x_{\bar{s}_{n_2}}) & 1\\
\vdots & & \vdots & \vdots\\
x_{\bar s_{1}}^{n_2-1}w_{1,2}(x_{\bar s_1}) & \ldots & x_{\bar
s_{n_2}}^{n_2-1}w_{1,2}(x_{\bar
s_{n_2}}) & y^{n_2-1}\\
\frac{w_{1,2}(x_{\bar s_{1}})}{z-x_{\bar s_1}} & \ldots & \frac{w_{1,2}(x_{\bar
s_{n_2}})}{z-x_{\bar s_{n_2}}} & \frac{1}{z-y}\end{vmatrix}
\end{multline}
where the sum is over all $\binom{n}{n_1}$ subsets $S\subset\{1,\ldots,n\}$
with $|S| = n_1$, where we write the elements of $S$ in increasing order:
$s_1<\ldots<s_{n_1}$, similarly for those of the complement $\bar{S} =
\{1,\ldots,n\}\setminus S$, $\bar{s}_1<\ldots<\bar{s}_{n_2}$, and where we
denote by $(-1)^{S}$ the sign of the permutation $s_1,\ldots,s_{n_1},\bar
s_1,\ldots,\bar s_{n_2}$. Note that the factor $(-1)^{n_2}$ of the previous
formula has disappeared.

Up to a diagonal scaling of the columns, both determinants in
\eqref{appendix:Lagrange} are Cauchy-Vandermonde determinants. By virtue of
\eqref{cauchy:vdm:det}, this yields
\begin{multline}\label{appendix:Lagrange2}
= \frac{1}{(z-y)^2}\left(\prod_{j=1}^n \frac{y-x_j}{z-x_j}\right) \sum_{S}
(-1)^{S}
\begin{vmatrix} w_{1,1}(x_{s_1}) & \ldots & w_{1,1}(x_{s_{n_1}}) \\
\vdots & & \vdots \\
x_{s_{1}}^{n_1-1}w_{1,1}(x_{s_1}) & \ldots &
x_{s_{n_1}}^{n_1-1}w_{1,1}(x_{s_{n_1}})
\end{vmatrix}\\
\begin{vmatrix} w_{1,2}(x_{\bar{s}_1}) & \ldots & w_{1,2}(x_{\bar{s}_{n_2}})\\
\vdots & & \vdots\\
x_{\bar s_{1}}^{n_2-1}w_{1,2}(x_{\bar s_1}) & \ldots & x_{\bar
s_{n_2}}^{n_2-1}w_{1,2}(x_{\bar s_{n_2}})\end{vmatrix}.
\end{multline}
But this is precisely what one gets from Lagrange expansion of the determinant
in \eqref{proof3:righthandside:reduced}.

\subsection*{Proof that \eqref{proof4:beforeclaim} implies
\eqref{proof4:righthandside:reduced}}

The proof that \eqref{proof4:beforeclaim} implies
\eqref{proof4:righthandside:reduced} follows in the same way. The only
difference is that there is an extra row and column in each Cauchy-Vandermonde
determinant, thereby yielding the pre-factor \eqref{proof4:prefactor}.


\begin{thebibliography}{99}


\bibitem{AvMV}
    M. Adler, P. van Moerbeke, and P. Vanhaecke,
    Moment matrices and multi-component KP, with applications to
    random matrix theory,
    Comm. Math. Phys. 286 (2009), 1--38.
\bibitem{Apt}
    A.I. Aptekarev,
    Multiple orthogonal polynomials,
    J. Comp. Appl. Math. 99 (1998), 423--447.
\bibitem{Baik}
    J. Baik,
    On the Christoffel-Darboux kernel for random Hermitian matrices with external
    source,
    manuscript 2009 (arXiv:0809.3970).
\bibitem{BDS}
    J. Baik, P. Deift and E. Strahov,
    Products and ratios of characteristic polynomials of random hermitian matrices,
    J. Math. Phys. 44 (2003) 3657--3670.
\bibitem{BK1}
    P.M. Bleher and A.B.J. Kuijlaars, Random matrices with
    external source and multiple orthogonal polynomials, Int.
    Math. Research Notices 2004, no 3 (2004), 109--129.
\bibitem{Bor}
    A. Borodin,
    Biorthogonal ensembles,
    Nucl. Phys. B536 (1999), 704--732.
\bibitem{BS}
    A. Borodin and E. Strahov,
    Averages of characteristic polynomials in random matrix theory,
    Comm. Pure Appl. Math. 59 (2006), 161--253.
\bibitem{BH3}
    E. Br\'ezin and S. Hikami, Characteristic polynomials of random matrices,
    Comm. Math. Phys. 214 (2000), 111--135.
    4140--4149.
\bibitem{BH1}
    E. Br\'ezin and S. Hikami, Universal singularity at the closure
    of the gap in a random matrix theory, Phys. Rev. E 57 (1998),
    4140--4149.
\bibitem{BH2}
    E. Br\'ezin and S. Hikami, Level spacing of random matrices in
    an external source, Phys. Rev. E 58 (1998), 7176--7185.
\bibitem{DK2}
    E. Daems and A.B.J. Kuijlaars,
    Multiple orthogonal polynomials of mixed type and non-intersecting Brownian
    motions, J. Approx. Theory 146 (2007), 91--114.
\bibitem{Dei}
    P. Deift, Orthogonal Polynomials and Random Matrices: a Riemann-Hilbert
    Approach. Courant Lecture Notes in Mathematics Vol. 3, Amer. Math. Soc.,
    Providence R.I. 1999.
\bibitem{Desrosiers1}
    P. Desrosiers and P.J. Forrester,
    A note on biorthogonal ensembles,
    J. Approx. Theory 152 (2008), 167--187.
\bibitem{Fid}
    U. Fidalgo Prieto, A. L\'opez Garc\'ia, G. L\'opez Lagomasino and V.N.
    Sorokin,
    Mixed type multiple orthogonal polynomials for two Nikishin systems,
    manuscript 2008 (arXiv:0812.1219).
\bibitem{FIK}
    A.S. Fokas, A.R. Its, and A.V. Kitaev,
    The isomonodromy approach to matrix models in 2D quantum gravity,
    Commun. Math. Phys. 147 (1992), 395--430.
\bibitem{FS}
    Y.V. Fyodorov and E. Strahov,
    An exact formula for general spectral correlation function of random
    Hermitian matrices,
    J. Phys. A: Math. Gen. 36 (2003), 3203--3213.
\bibitem{Gant}
    F.R. Gantmacher,
    The Theory of Matrices, Vol. 1,
    Chelsea Publishing Company, New York, 1959.
\bibitem{Gelfand}
    I. Gelfand, S. Gelfand, V. Retakh and R. Wilson,
    Quasideterminants,
    Advances in Mathematics 193 (2005), 56-–141.
\bibitem{GVL}
    G.H. Golub and C.F. Van Loan,
    Matrix Computations, The Johns
    Hopkins University Press, third edition, 1996.
\bibitem{KMcG}
    S. Karlin and J. McGregor, Coincidence probabilities,
    Pacific J. Math., 9 (1959), 1141--1164.
\bibitem{KS}
    J.P. Keating and N.C. Snaith,
    Random Matrix Theory and $\zeta(1/2 + it)$,
    Commun. Math. Phys. 214 (2000), 57–-89.
\bibitem{Kui}
    A.B.J. Kuijlaars,
    Multiple orthogonal polynomial ensembles,
    manuscript 2009 (arXiv:0902.1058).
\bibitem{Miranian}
    L. Miranian,
    Matrix-valued orthogonal polynomials on the real line: some extensions of
    the classical theory,
    J. Phys. A: Math. Gen. 38 (2005), 5731--5749.
\bibitem{SVI}
    V.N. Sorokin and J. Van Iseghem,
    Algebraic aspects of matrix orthogonality for vector polynomials,
    J. Approx. Theory 90 (1997), 97--116.
\bibitem{SF}
    E. Strahov, Y.V. Fyodorov,
    Universal results for correlations of characteristic polynomials: Riemann-Hilbert
    approach,
    Comm. Math. Phys. 241 (2003), 343--382.
\bibitem{Sz}
     G. Szeg\H{o},
     Orthogonal polynomials,
     Amer. Math. Soc. Coll. Publ. Vol 23, Amer. Math. Soc., Providence, R.I. 1975.
\bibitem{Uv}
    V. B. Uvarov,
    The connection between systems of polynomials orthogonal with respect to
    different distribution functions,
    USSR Comput. Math. and Math. Phys. 9 (Part 2) (1969), 25--36.
\bibitem{VA}
    W. Van Assche,
    Multiple orthogonal polynomials, irrationality and transcendence,
    in ``Continued fractions: from analytic number theory to constructive
    approximation'',
    Contemporary Mathematics 236 (1999), 325--342.
\bibitem{VAGK} W. Van Assche, J.S. Geronimo, and A.B.J. Kuijlaars,
    Riemann-Hilbert problems for multiple orthogonal polynomials,
     Special Functions 2000: Current Perspectives and Future Directions
    (J. Bustoz et al., eds.), Kluwer, Dordrecht, 2001, pp. 23--59.
\bibitem{Zinn}
    P. Zinn-Justin,
    Universality of correlation functions of Hermitian random matrices in an
    external field,
    Comm. Math. Phys. 194 (1998), 631--650.
\end{thebibliography}
\end{document}